\documentclass[12pt]{article}
\usepackage{amsmath}
\usepackage{amsfonts}
\usepackage{latexsym, amssymb, amsmath, amscd, amsfonts, epsfig, graphicx,texdraw}
\usepackage{mathrsfs}
\usepackage{colordvi}
\usepackage{ifpdf}
\usepackage{lscape}
\newtheorem{theorem}{Theorem}[section]

\newtheorem{lemma}[theorem]{Lemma}
\newtheorem{corollary}[theorem]{Corollary}
\newtheorem{prop}[theorem]{Proposition}
\newtheorem{definition}[theorem]{Definition}

 \numberwithin{equation}{section}

\def\qed{\hfill \rule{4pt}{7pt}}

\def\pf{\noindent {\it Proof.} }

\begin{document}
\begin{center}
{\large {\bf Equivalence Classes of Full-Dimensional

$0/1$-Polytopes with Many Vertices}}

\vskip 6mm {\small William Y.C. Chen$^1$ and Peter L. Guo$^2$\\[%
2mm] Center for Combinatorics, LPMC-TJKLC\\
Nankai University, Tianjin 300071,
P.R. China \\[3mm]
$^1$chen@nankai.edu.cn, $^2$lguo@cfc.nankai.edu.cn \\[0pt%
] }
\end{center}
\begin{abstract}
Let $Q_n$ denote the $n$-dimensional hypercube with the vertex set
$V_n=\{0,1\}^n$. A $0/1$-polytope of  $Q_n$ is a convex hull of a
subset of $V_n$. This paper is concerned with the enumeration of
equivalence classes of full-dimensional  $0/1$-polytopes under the
symmetries of the hypercube.  With the aid of a computer program,
Aichholzer completed the enumeration of equivalence classes of
full-dimensional $0/1$-polytopes for $Q_4$, $Q_5$, and those of
$Q_6$ up to $12$ vertices. In this paper, we present a method to
compute the number of equivalence classes of full-dimensional
0/1-polytopes of $Q_n$ with more than $2^{n-3}$ vertices. As an
application, we finish the counting of equivalence classes of
full-dimensional $0/1$-polytopes of $Q_6$ with more than 12
vertices.
\end{abstract}
\vskip 3mm



\allowdisplaybreaks

\section{Introduction}

Let $Q_n$ denote the  $n$-dimensional hypercube
with vertex set $V_n=\{0,1\}^n$.
A $0/1$-polytope of $Q_n$ is defined
to be the convex hull of a subset of $V_n$.
 The study of $0/1$-polytopes has drawn much attention
 from different points of view, see, for example,  \cite{BP,Bil,Fle,GK,Haiman,Kor,Zon2},
see also the survey of Ziegler \cite{Zei4}.

In this paper, we are concerned with
the problem of  determining
the number of equivalence classes
of $n$-dimensional $0/1$-polytopes of $Q_n$
 under the  symmetries of $Q_n$, which has been
 considered as a difficult problem, see Ziegler \cite{Zei4}.
  It is also listed by
Zong \cite[Problem 5.1]{Zon2} as  one of the fundamental problems
concerning $0/1$-polytopes.

An $n$-dimensional $0/1$-polytope of $Q_n$ is also called
a full-dimensional $0/1$-polytope of $Q_n$.
Two  $0/1$-polytopes are said to be
equivalent if one can be transformed to the other by a symmetry
of $Q_n$. Such a equivalence relation is also called  the
$0/1$-equivalence relation. Figure \ref{fig} gives
representatives of $0/1$-equivalence classes of $Q_2$, among
which (d) and (e) are
full-dimensional.

\begin{figure}
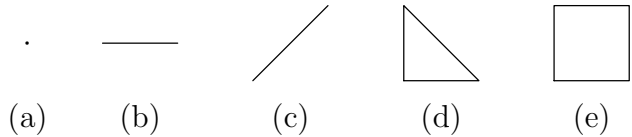

\centertexdraw{
\drawdim mm \linewd 0.15
\setgray 0

\move(10 15) \lvec(20 15)
\move(50 10)\lvec(60 10)\lvec(50 20)\lvec(50 10)
\move(70 10)\lvec(80 10)\lvec(80 20)\lvec(70 20)\lvec(70 10)
\textref h:C v:C \htext(0 15){$\cdot$}
\move(30 10)\lvec(40 20)
\textref h:C v:C \htext(0 5){(a)}\textref h:C v:C \htext(15 5){(b)}
\textref h:C v:C \htext(35 5){(c)}\textref h:C v:C \htext(55 5){(d)}
\textref h:C v:C \htext(75 5){(e)}
}
\caption{0/1-Polytopes of the square}
\end{figure}\label{fig}

Sarangarajan and Ziegler \cite[Proposition 8]{Zei4}
found an lower bound on the number of equivalence
classes of full-dimensional $0/1$-polytopes of $Q_n$.
As far as exact enumeration is concerned,
full-dimensional
0/1-equivalence classes of $Q_4$ were counted by
Alexx Below,
see Ziegler \cite{Zei4}. With the aid of a computer
program, Aichholzer \cite{Aic2}
 completed the
enumeration  of full-dimensional $0/1$-equivalence classes of
of $Q_5$, and  those of
$Q_6$ up to $12$ vertices, see Aichholzer \cite{Aic3} and
Ziegler \cite{Zei4}. The 5-dimensional hypercube $Q_5$
has been considered as the last case that one can hope for
a complete solution to the enumeration of full-dimensional
$0/1$-equivalence classes.

The objective of this paper is to present a method to compute the number of
full-dimensional $0/1$-equivalence classes of $Q_n$ with more than $2^{n-3}$ vertices. As an application, we
solve the
enumeration problem for full-dimensional $0/1$-equivalence classes of
the 6-dimensional hypercube with more than 12 vertices.

To describe our approach, we introduce some notation.
 Denote by
$\mathcal{A}_n(k)$ (resp., $\mathcal{F}_n(k)$) the set of  (resp.,
full-dimensional) $0/1$-equivalence classes of $Q_n$ with $k$
vertices. Let $\mathcal{H}_n(k)$ be the
set of  $0/1$-equivalence classes of  $Q_n$ with
$k$ vertices that are not full-dimensional.
The cardinalities of $\mathcal{A}_n(k)$,
$\mathcal{F}_n(k)$ and $\mathcal{H}_n(k)$
are denoted respectively
by $A_n(k)$, $F_n(k)$ and $H_n(k)$.
It is clear that  any
full-dimensional $0/1$-polytope of $Q_n$ has at least $n+1$
vertices, i.e., $F_n(k)=0$ for $1\leq k\leq n$.

The  starting point of this paper is the following
obvious relation
\begin{equation}\label{cg10}
F_n(k)=A_n(k)-H_n(k).
\end{equation}
The number $A_n(k)$
can be computed based on the cycle index of
the hyperoctahedral
group.  We can deduce that  $H_n(k)=0$ for $k>2^{n-1}$ based on a result duo to Saks. For the purpose of computing  $H_n(k)$ for $2^{n-2}<k\leq 2^{n-1}$, we transform the computation of  $H_n(k)$
to the determination of the number of
equivalence classes of $0/1$-polytopes
with $k$ vertices
that are contained in the spanned hyperplanes of $Q_n$.
To be more specific, we show that $\mathcal{H}_n(k)$ for
$2^{n-2}<k\leq 2^{n-1}$ can be decomposed
into a disjoint union of equivalence classes of
$0/1$-polytopes
that are contained  in the spanned hyperplanes of $Q_n$.
In particular,
 for $n=6$ and
$k>16$, we obtain  the number of
 full-dimensional $0/1$-equivalence classes
of $Q_6$ with $k$ vertices.

Using a similar idea as in the case $2^{n-2}<k\leq 2^{n-1}$,
we can compute  $H_n(k)$ for   $2^{n-3}<k\leq 2^{n-2}$.
For $n=6$ and $13\leq k \leq 16$,
we obtain the number of
 full-dimensional $0/1$-equivalence classes
of $Q_6$ with $k$ vertices.
 Together with the computation of
Aichholzer up to 12 vertices,
we have completed the enumeration of
 full-dimensional $0/1$-equivalence classes  of the
6-dimensional hypercube.

\section{The cycle index of the hyperoctahedral group}

The group of symmetries of $Q_n$ is known as
the hyperoctahedral
group $B_n$. In this section,
we review the cycle index of $B_n$
acting on the vertex set $V_n$.
Since $0/1$-equivalence classes  of $Q_n$ coincide with  nonisomorphic vertex
colorings of $Q_n$ by
 using two colors, we may compute the number $A_n(k)$ from
 the cycle index of $B_n$.

Let $G$ be a   group acting on a
finite set $X$.
For any $g\in G$, $g$ induces a permutation on $X$.
The cycle type of a permutation is defined to be a multiset $\{1^{c_1},2^{c_2},\ldots\}$, where $c_i$ is the number of cycles of length $i$ that appear in the cycle decomposition of the permutation.
For $g\in G$, denote
by $c(g)=\{1^{c_1},2^{c_2},\ldots\}$ the cycle type of the
permutation on $X$ induced by $g$.  Let
$z=(z_1,z_2,\ldots)$ be a sequence of indeterminants, and let
\[
z^{c(g)}=z_1^{c_1}z_2^{c_2}\cdots.\] The cycle
index of $G$ is defined as follows
\begin{equation}\label{cg22}
Z_G(z)=Z_G(z_1,z_2,\ldots)=\frac{1}{|G|}\sum_{g\in
G}z^{c(g)}.
\end{equation}
According to  P\'olya's theorem, the cycle
index in (\ref{cg22})
can be applied to count nonisomorphic  colorings of
$X$ by using a given number of colors.

For a vertex coloring of $Q_n$ with two colors, say,
 black and white,  the black vertices can be considered as
 vertices of a $0/1$-polytope of $Q_n$. This establishes
a one-to-one correspondence between  equivalence classes of vertex
colorings and  $0/1$-equivalence classes  of $Q_n$.
Let $Z_n(z)$ denote
the cycle index of $B_n$ acting on the vertex set $V_n$.
Then, by P\'olya's theorem
\begin{equation}\label{cg13}
A_n(k)=\left[u_1^ku_2^{2^n-k}\right]C_n(u_1,u_2),
\end{equation}
 where $C_n(u_1,u_2)$ is the polynomial obtained from $Z_n(z)$
 by substituting    $z_i$ with $u_1^i+u_2^i$, and
 $[u_1^pu_2^q]C_n(u_1,u_2)$ denotes the coefficient of
 $u_1^pu_2^q$ in $C_n(u_1,u_2)$.

Clearly, the total number of $0/1$-equivalence classes
of $Q_n$ is given by
\begin{equation}\label{cg12}
\sum_{k=1}^{2^n}A_n(k)=C_n(1,1).
\end{equation} It should be noted that $C_n(1,1)$ also
equals  the number of
types of Boolean functions, see Chen \cite{Che1} and
references therein. This number is also related to
 configurations
of $n$-dimensional
Orthogonal Pseudo-Polytopes, see, e.g.,
 Aguila \cite{Agu}.
The computation of $Z_n(z)$ has been
studied by  Chen \cite{Che1},  Harrison
and High \cite{Har}, and P\'olya \cite{Polya}, etc.
 Explicit expressions of $Z_n(z)$
for $n\leq 6$
can be found in \cite{Agu}, and we list them bellow.
\begin{align*}
Z_1(z)&=z_1,\\[5pt]
Z_2(z)&=\frac{1}{8}\left(\begin{array}{l}
z_1^4+2z_1^2z_2+3z_2^2+2z_4
\end{array}\right),\\[5pt]
Z_3(z)&=\frac{1}{48}\left(\begin{array}{l}
z_1^8+6z_1^4z_2^2+13z_2^4+8z_1^2z_3^2+12z_4^2+8z_2z_6
\end{array}\right),\\[5pt]
Z_4(z)&=\frac{1}{384}\left(\begin{array}{l}
z_1^{16}+12z_1^8z_2^4+12z_1^4z_2^6+51z_2^8+48z_8^2\\[5pt]
+48z_1^2z_2z_4^3+84z_4^4+96z_2^2z_6^2+32z_1^4z_3^4
\end{array}\right),\\[9pt]
Z_5(z)&=\frac{1}{3840}\left(\begin{array}{l}
 z_1^{32}+20z_1^{16}z_2^8+ 60z_1^8z_2^{12}+231z_2^{16}+
80z_1^8z_3^8+240z_1^4z_2^2z_4^6\\[5pt]
 +240z_2^4z_4^6+520z_4^8+384z_1^2z_5^6+160z_1^4z_2^2z_3^4z_6^2+
720z_2^4z_6^4\\[5pt]
 +480z_8^4+384z_2z_{10}^3+320z_4^2z_{12}^2
\end{array}\right),\\[9pt]
Z_6(z)&=\frac{1}{46080}\left(\begin{array}{l}
z_1^{64}+30z_1^{32}z_2^{16}+180z_1^{16}z_2^{24}+120z_1^8z_2^{28}+
1053z_2^{32}+160z_1^{16}z_3^{16}+\\[5pt]
640z_1^4z_3^{20}+720z_1^8z_2^4z_4^{12}+
1440z_1^4z_2^6z_4^{12}+2160z_2^8z_4^{12}+4920z_4^{16}+\\[5pt]
2304z_1^4z_5^{12}+
960z_1^8z_2^4z_3^8z_6^4+5280z_2^8z_6^8+3840z_1^2z_2z_3^2z_6^9+5760z_8^8\\[5pt]
+1920z_2^2z_6^{10}+6912z_2^2z_{10}^6+3840z_4^4z_{12}^4+
3840z_4z_{12}^5
\end{array}\right).
\end{align*}

The method  of Chen for
 computing $Z_n(z)$ is
 based on the cycle structure of a power of a signed permutation.
Let us recall the notation of a signed permutation.
A signed permutation on $\{1,2,\ldots,n\}$ is
a permutation on $\{1,2,\ldots,n\}$
with a $+$ or a $-$ sign
attached to each element $1,2,\ldots,n$.
Following the notation
in Chen \cite{Che1} or  Chen and Stanley
\cite{Che2}, we may write a
signed permutation in terms of the cycle decomposition
and ignore the plus sign $+$.
For example, $(\overline{2}4\overline{5})(3)(1\overline{6})$
represents
a signed permutation, where $(245)(3)(16)$ is
called  its underlying permutation.
The action of  a signed permutation $w$ on the vertices of $Q_n$
is
defined as follows.  For a vertex $(x_1, x_2,\ldots, x_n)$ of $Q_n$,
we define
$w((x_1,x_2,\ldots,x_n))$ to be the vertex  $(y_1,y_2,\ldots,y_n)$
as given by
\begin{equation}\label{cg19}
y_i=\left\{\begin{array}{ll}
x_{\pi(i)}, &\mbox{if $i$ has the sign $+$},\\[5pt]
1-x_{\pi(i)}, &\mbox{if $i$ has the sign $-$},\end{array}\right.
\end{equation} where $\pi$ is the  underlying
permutation of $w$.

For the purpose of this paper, we define the cycle
type of a signed permutation  $w\in B_n$ as the cycle
type of its
underlying permutation. For example,
$(\overline{2}4\overline{5})(3)(1\overline{6})(7)$
has cycle type $\{1^2,2,3\}$. We should note
that the above definition of a cycle type of a signed permutation
is different from the definition in terms of double
partitions as in \cite{Che1} because it will be shown in Section 5 that
any  signed permutation
that fixes a  spanned hyperplane of $Q_n$ either have all positive cycles or
all negative cycles.

We end this section with the following
formula of  Chen
\cite{Che1}, which will be used in
 Section 6 to compute the cycle index of the group that fixes  a spanned hyperplane of $Q_n$.

\begin{theorem}\label{cg7}
Let $G$ be a group that acts on some finite
set $X$. For any $g\in G$, the number
of $i$-cycles of the  permutation on $X$
induced by $g$  is given by
\[\frac{1}{i}\sum_{j|i}\mu(i/j)\psi(g^j),\]
where $\mu$ is the
classical number-theoretic M\"{o}bius
function and $\psi(g^j)$ is
the number of fixed points of $g^j$ on $X$.
\end{theorem}

\section{$0/1$-Polytopes with many vertices}

In this section, we find an inequality concerning the dimension of a
$0/1$-polytope of $Q_n$ and the  number of its vertices.
This inequality plays a key role in the computation
 of $F_n(k)$ for $k>2^{n-3}$.

The main theorem of this section is given below.

\begin{theorem}\label{cg14}
Let $P$ be a $0/1$-polytope of $Q_n$ with more than $2^{n-s}$ vertices,
where $1\leq s\leq n$. Then  we have \[\mathrm{dim}(P)\geq n-s+1.\]
\end{theorem}

The above theorem can be  deduced from the following assertion.

\begin{theorem} \label{cgx}
For any $1\leq s \leq n$, the intersection of $s$ hyperplanes
in $\mathbb{R}^n$ with linearly independent
normal vectors contains at
most $2^{n-s}$ vertices of $Q_n$.
\end{theorem}

Indeed, it is not difficult to see that Theorem \ref{cgx}
implies Theorem \ref{cg14}. Let $P$ be a $0/1$-polytope of $Q_n$ with more than $2^{n-s}$ vertices.
Suppose to the contrary that $\mathrm{dim}(P)\leq n-s$.
It is known that
the affine space spanned by $P$ can be expressed as the intersection of a collection of  hyperplanes. Since $\mathrm{dim}(P)\leq n-s$,
there  exist $s$
hyperplanes $H_1,H_2,\ldots,H_s$ whose normal vectors are linearly independent such
that the intersection of $H_1, H_2, \ldots, H_s$ contains $P$.
 Let $V(P)$ denote the vertex set of $P$. By Theorem \ref{cgx}, we have
 \[|V(P)|\leq \left|\left(\bigcap_{i=1}^s
H_i\right)\bigcap V_n\right|\leq 2^{n-s},\] which is a contradiction to
 the assumption   that $P$ contains more than $2^{n-s}$ vertices of $Q_n$. So we conclude that $\mathrm{dim}(P)\geq n-s+1$.

\noindent\textit{Proof of Theorem \ref{cgx}.}
Assume that, for
$1\leq i\leq s$, \[H_i\colon a_{i1}x_1+a_{i2}x_2+
\cdots +a_{in}x_n=b_i\] are $s$ hyperplanes
in $\mathbb{R}^n$, whose
normal vectors ${a}_i=(a_{i1},
\ldots,a_{in})$ are linearly independent.
We aim to show that the intersection of $H_1,H_2,\ldots,H_s$ contains at
most $2^{n-s}$ vertices of $Q_n$.
We may express the intersection of  $H_1,H_2,\ldots,H_s$ as the solution
of a system of linear equations, that is,
\begin{equation}\label{cg50x}
\bigcap_{i=1}^s H_i=\{{x^T}\colon
A{x}=b\},\end{equation}
where  $A$ denotes the matrix $(a_{ij})_{1\leq i\leq s,1\leq j \leq n}$, $x=(x_1, x_2, \ldots, x_n)^T$, and ${b}=(b_1,\ldots,b_s)^T$, $T$ denotes the transpose of a
vector.
Then Theorem \ref{cgx} is equivalent to the following inequality
\begin{equation}\label{cg50}
\left|V_n\bigcap \left\{{x^T}\colon
A{x}={b}\right\}
\right|\leq 2^{n-s}.\end{equation}

We now proceed to prove  (\ref{cg50}) by induction
on $n$ and $s$. We first consider the case  $s=1$. Suppose that  $H\colon c_1x_1+c_2x_2+\cdots+c_nx_n=c$ is a hyperplane in
$\mathbb{R}^n$.  Assume that among the coefficients $c_1, c_2, \ldots, c_n$ there are  $i$ of them
that are nonzero.  Without loss of generality, we may assume that $c_1, c_2,\ldots, c_{i}$ are nonzero, and $c_{i+1}= c_{i+2}=\cdots=
c_{n}=0$.  Clearly, $H$ reduces to a hyperplane in the $i$-dimensional Euclidean space $\mathbb{R}^{i}$. Such a hyperplane with nonzero coefficients is called a skew hyperplane.
Now the vertices of $Q_n$ contained in
$H$ are of the form $(d_1,\ldots,d_{i},d_{i+1},\ldots,d_n)$ where $(d_1,\ldots,d_{i})$ are vertices of $Q_{i}$ contained in the
skew hyperplane $H'\colon c_1x_1+c_2x_2+\cdots c_ix_i=b$.
Clearly, for each vertex $(d_1,d_2,\ldots, d_i)$ in $H'$, there are
 $2^{n-i}$ choices for $(d_{i+1},d_{i+2},\ldots, d_n)$
such that $(d_1,d_2,\ldots, d_n)$ is contained in $H$.
Using Sperner's lemma (see,
for example, Lubell \cite{ Lub}), Saks \cite[Theorem 3.64]{Sak} has
shown  that
the  number of vertices of $Q_i$ contained in
a skew hyperplane
does not exceed ${i \choose \lfloor \frac{i}{2}\rfloor}$. Let \[ f(n,i)=2^{n-i}{i \choose \lfloor \frac{i}{2}\rfloor}.\]
Thus the number of vertices of $Q_n$ contained in  $H$ is at most
$f(n,i)$. It is easy to check  that
\[  {f(n,i) \over f(n,i+1)}=
\left\{\begin{array}{ll}
\frac{i+2}{i+1}, &\mbox{if $i$ is even},\\[5pt]
1, &\mbox{if $i$ is odd}.\end{array}\right
.\] This yields  $f(n,i)\geq f(n,i +1)$ for any $i=1,2,\ldots, n-1$.
Hence $H$ contains at most $f(n,1)=2^{n-1}$ vertices of $Q_n$, which implies
 (\ref{cg50}) for $s=1$.

We now consider the case $s=n$. In this case, since  the normal vectors $a_1,\ldots,a_n$ are linearly independent, the square matrix $A$ is nonsingular. It follows that
$A{x}={b}$ has exactly one solution. Therefore,  inequality (\ref{cg50}) holds when $s=n$.

So we are left with cases of
$n,s$ such that $1<s<n$. We shall use induction to complete the proof.
Suppose that
(\ref{cg50}) holds for $n', s'$ such that $n'\leq n, s'\leq s$ and
$(n',s')\neq (n,s)$.

Since the normal vector $a_1$ is nonzero, there exists some
$j_0$ ($1\leq j_0\leq n$) such that $a_{1j_0}\neq 0$. Without loss of generality, we may
assume  $a_{ij_0}=0$ for $2\leq i\leq
s$ since one can apply elementary row
transformations to the system of linear
equations $A{x}={b}$ to ensure that the assumption is valid. For a vector
${v}=(v_1,\ldots,v_n)\in \mathbb{R}^n$,
let
${v}^j=(v_1,\ldots,v_{j-1},v_{j+1},
\ldots,v_{n})\in \mathbb{R}^{n-1}$ be the
vector obtained from ${v}$ by
deleting the $j$-th coordinate.
We now have two cases.

Case 1. The vectors ${a}_1^{j_0},\ldots,
{a}_s^{j_0}$ are linearly
dependent. Since ${a}_2,\ldots,
{a}_s$ are linearly
independent and $a_{ij_0}=0$ for $2\leq i\leq s$,
it is clear that
${a}_2^{j_0},\ldots,{a}_s^{j_0}$ are
linearly independent.
So the vector $a_1^{j_0}$ can be expressed as a linear combination
of ${a}_2^{j_0},\ldots,
{a}_s^{j_0}$. Assume that ${a}_1^{j_0}=\alpha_2{a}_2^{j_0}+
\ldots+\alpha_s{a}_s^{j_0}$, where $\alpha_k\in \mathbb{R}$ for
$2\leq k\leq s$.
For $2\leq k\leq s$,  multiplying the $k$-th row  by $\alpha_k$ and subtracting
 it from the first row, then the first equation $a_{11}x_1+a_{12}x_2+\cdots+a_{1n}x_n=b_1$ becomes
  \begin{equation}\label{cg52}
 a_{1j_0}x_{j_0}=b_1-\sum_{k=2}^s\alpha_kb_k.\end{equation}
 Let
$A'=({a}_2,\ldots,{a}_s)^T$ and $b'=(b_2,\ldots,b_s)^T$.
Note that all
 entries in  the $j_0$-th column  of $A'$ are  zero since we have assumed  $a_{ij_0}=0$ for $2\leq i\leq
s$. Let
$A'_{j_0}$
be the matrix
obtained from $A'$ by removing this zero column.
From equation (\ref{cg52}),  the value $x_{j_0}$ in the $j_0$-th coordinate of the solutions of $A{x}={b}$
is  fixed.
Then solutions of $ A{x}={b}$ can be obtained from the solutions of
$A'_{j_0}x=b'$ by adding the value of $x_{j_0}$ to  the $j_0$-th coordinate.
 Concerning the number of vertices of $Q_n$ contained in $\left\{{x}^T\colon
A{x}={b}\right\}$, we consider the following two cases.

(1). The  value $x_{j_0}$ is not equal to $0$ or $1$.
In this case, no vertex of $Q_n$ is contained in $\left\{{x}^T\colon
A{x}={b}\right\}$. Hence  inequality (\ref{cg50}) holds.

(2). The value $x_{j_0}$ is equal to $0$ or $1$. Since every vertex of $Q_n$ contained in $\left\{{x^T}\colon
A{x}={b}\right\}$ is obtained from a vertex of $Q_{n-1}$ contained in $
\left\{x^T\colon
A'_{j_0}x={b'}
\right\}$ by adding $x_{j_0}$ in the $j_0$-th coordinate,
it follows that
\begin{equation}\label{ab}
\left|V_n\bigcap \left\{{x^T}\colon
A{x}={b}\right\}
\right|=\left|V_{n-1}\bigcap\,
\left\{{x^T}\colon
A'_{j_0}{x}={b'}
\right\}\right|.\end{equation}
By the induction hypothesis, we find
\[
\left|V_{n-1}\bigcap\,
\left\{{x^T}\colon
A'_{j_0}{x}=b'
\right\}\right|\leq
2^{(n-1)-(s-1)}=2^{n-s}.
\]
In view of (\ref{ab}), we obtain (\ref{cg50}).

Case 2. Suppose   ${a}_1^{j_0},\ldots,
{a}_s^{j_0}$ are linearly
independent. Assume that  the value of $x_{j_0}$ in the solutions of $\left\{{x^T}\colon
A{x}=b\right\} $ can be taken
$0$ or $1$. Then the vertices of $Q_n$ contained in $\left\{{x^T}\colon
A{x}=b\right\} $
can be decomposed into
a disjoint union of the following two sets
\[S_0=V_n\bigcap\left\{{x^T}\colon x_{j_0}=0,\
A{x}=b\right\}\]and
\[S_1=V_n\bigcap\left\{{x^T}\colon x_{j_0}=1,\
A{x}={b}\right\}.\]
We first consider the set $S_0$. Let $A''=({a}_1^{j_0},\ldots,
{a}_s^{j_0})^T$
be the matrix obtained from $A$
by deleting the ${j_0}$-th column. Then vertices of $Q_n$ contained in $S_0$ are obtained from  the vertices of $Q_{n-1}$ contained in $\{x^T\colon A''x={b}\}$ by adding $0$ to the $j_0$-th coordinate. So we have
\[|S_0|=\left|
V_{n-1}\bigcap\left\{{x^T}\colon \
A''x={b}\right\}\right|\leq 2^{n-1-s},\] where the inequality follows from the induction
hypothesis. Similarly,
we   get $|S_1|\leq
2^{n-1-s}$. Hence \[\left|V_n\bigcap \left\{{x^T}\colon
A{x}={b}\right\}\right|=|S_0|+|S_1|\leq
2^{n-s}.\]

Combining the above two cases,  inequality (\ref{cg50}) is true for
 $1\leq s\leq n$.
This completes the proof. \qed

Note that the upper bound $2^{n-s}$ is sharp. For example, it is
easy to see the intersection of hyperplanes $x_i=0$ ($1\leq i\leq
s$) contains exactly $2^{n-s}$ vertices of $Q_n$.

By Theorem \ref{cgx}, we see that
every $0/1$-polytope of $Q_n$ with more than $2^{n-1}$ vertices is full-dimensional. As a direct consequence, we obtain the following relation.

\begin{corollary}\label{cg51}
For $k>2^{n-1}$, we have
\[F_n(k)=A_n(k).\]
\end{corollary}

Form Corollary \ref{cg51}, the number $F_n(k)$ for $k>2^{n-1}$ can be computed from the cycle index of the hyperoctahedral group, that is, for $k>2^{n-1}$
\begin{equation*}
F_n(k)=\big[u_1^ku_2^{2^n-k}\big]C_n(u_1,u_2).
\end{equation*}
For $n=4,5$ and $6$, the values of
$F_n(k)$ for $k>2^{n-1}$ are given
in Tables \ref{tt1}, \ref{tt3} and \ref{tt5}.
\begin{table}[h,t]
\begin{center}
\begin{tabular}{|l|l|l|l|l|l|l|l|l|l|}\hline
$k$&9&10&11&12&13&14&15&16\\\hline
$F_4(k)$&56&50&27&19&6&4&1&1\\\hline
\end{tabular}
\end{center}
\caption{$F_4(k)$ for $k>8$.}\label{tt1}
\end{table}

\begin{table}[h,t]
\begin{center}
\begin{tabular}{|l|l|l|l|l|l|l|l|l|l|}
 \hline
$k$&17&18&19&20&21&22&23&24\\\hline
$F_5(k)$&158658&133576&98804&65664&38073&19963&9013&3779\\\hline
$k$&25&26&27&28&29&30&31&32\\\hline $F_5(k)$&1326&472&
131&47&29&5&1&1
\\\hline
\end{tabular}
\end{center}
\caption{$F_5(k)$ for $k>16$.}\label{tt3}
\end{table}

\begin{table}[h,t]
\begin{center}
\begin{tabular}{|p{10mm}|l|p{10mm}|l|p{10mm}|l|}\hline
$k$&$F_6(k)$&$k$&$F_6(k)$\\\hline
33&38580161986426&49&3492397119\\\hline34&35176482187398&
50&1052201890
\\\hline35&30151914536933&51&290751447\\\hline
36&24289841497881&52&73500514\\\hline37&18382330104696
&53&16938566\\\hline38&13061946976545&
54&3561696\\\hline39&8708686182967&55&681474
\\\hline40&5443544478011&56&120843\\\hline
41&3186944273554&57&19735\\\hline42&1745593733454&58&3253\\\hline
43&893346071377&59&497\\\hline44
&426539774378&
60&103
\\\hline45&189678764492&61&16\\\hline46&78409442414&62&6\\\hline
47&30064448972&63&1\\\hline48&10666911842&64&1\\\hline
\end{tabular}
\end{center}
\caption{$F_6(k)$ for $k>32$.}\label{tt5}
\end{table}

\section{$H_n(k)$ for $2^{n-2}<k\leq 2^{n-1}$}

In this section, we shall aim to compute
 $H_n(k)$ for $2^{n-2}<k\leq 2^{n-1}$.
We shall show that in this case the number $H_n(k)$
is determined by the number of (partial) $0/1$-equivalence
classes of a spanned hyperplane of $Q_n$ with $k$ vertices.  To this end, it is necessary to
consider all
possible spanned hyperplanes of $Q_n$. More precisely, we need
representatives of equivalence classes of such spanned hyperplanes.

Recall that a spanned hyperplane of $Q_n$ is a hyperplane
in $\mathbb{R}^n$ spanned by $n$ affinely
independent vertices of $Q_n$, that is, the affine space spanned by the vertices of $Q_n$ contained in this hyperplane is of dimension $n-1$.
 Let \[ H\colon a_1x_1+a_2x_2+\cdots+a_nx_n=b\] be a spanned
hyperplane of $Q_n$, where $|a_1|,\ldots,|a_n|, |b|$ are positive integers with greatest common divisor $1$. Let
\[\mathrm{coeff}(n)=\mathrm{max}\{|a_1|,\ldots,|a_n|\}.\]
It is clear that $\mathrm{coeff}(2)=\mathrm{coeff}(3)=1$.
The study of upper and lower bounds on the number $\mathrm{coeff}(n)$ has drawn
 much attention, see, for example, \cite{Aic1,AW,Hastad,Zei4}.
The following are known bounds on $\mathrm{coeff}(n)$ and $|b|$, see, e.g., \cite[Corollary 26]{Zei4} and \cite[Theorem 5]{Aic1},
\begin{equation*}
\frac{(n-1)^{(n-1)/2}}{2^{2n+o(n)}}\leq \mathrm{coeff}(n)\leq \frac{n^{n/2}}{2^{n-1}}\ \ \ \ \ \mbox{and}\ \ \ \ \
|b|\leq 2^{-n}(n+1)^{\frac{n+1}{2}}.
\end{equation*}

Using the above bounds,
Aichholzer  and Aurenhammer \cite{Aic1} obtained the exact values
of $\mathrm{coeff}(n)$ for $n\leq 8$
by computing all possible  spanned hyperplanes
of $Q_n$ up to dimension $8$. For example, they showed that
  $\mathrm{coeff}(4)=2$, $\mathrm{coeff}(5)=3$, and $\mathrm{coeff}(6)=5$.

As will be seen, in order to
compute $H_n(k)$ for $2^{n-2}<k\leq 2^{n-1}$, we need to
consider equivalence classes of
spanned hyperplanes of
$Q_n$ under the symmetries of $Q_n$.  Note that
the symmetries of $Q_n$ can be expressed by permuting the
coordinates and changing  $x_i$ to $1-x_i$ for some indices $i$.
Therefore, for each equivalence class of spanned hyperplanes
of $Q_n$, we can choose a representative
of the following form
\begin{equation}\label{hyper}
a_1x_1+a_2x_2+\cdots+a_tx_t=b,
\end{equation}
where $t\leq n $ and $0<a_1\leq a_2\leq \cdots\leq a_t\leq \mathrm{coeff}(n)$.

A complete list of spanned hyperplanes of $Q_n$ for $n\leq 6$ can been found in \cite{Aic}. The following hyperplanes are representatives of
equivalence classes of spanned hyperplanes of $Q_4$:
\begin{align*}
&x_1=0,\\
&x_1+x_2=1,\\
&x_1+x_2+x_3=1,\\
&x_1+x_2+x_3+x_4=\mbox{$1$ or $2$},\\
&x_1+x_2+x_3+2x_4=2.
\end{align*}
In addition to the above hyperplanes of $\mathbb{R}^4$, which
can also be viewed as spanned hyperplanes of $Q_5$,
we have the following representatives of equivalence classes
of spanned hyperplanes of  $Q_5$:
\begin{align*}
&x_1+x_2+x_3+x_4+x_5=\mbox{$1$ or $2$},\\
&x_1+x_2+x_3+x_4+2x_5=\mbox{$2$ or $3$},\\
&x_1+x_2+x_3+2x_4+2x_5=\mbox{$2$ or $3$},\\
&x_1+x_2+2x_3+2x_4+2x_5=\mbox{$3$ or $4$},\\
&x_1+x_2+x_3+x_4+3x_5=3,\\
&x_1+x_2+x_3+2x_4+3x_5=3,\\
&x_1+x_2+2x_3+2x_4+3x_5=4.
\end{align*}

When $n=6$, for the purpose of computing $F_6(k)$
for $16<k\leq 32$, we need the representatives of equivalence
classes of spanned hyperplanes of $Q_6$ containing more
than $16$ vertices of $Q_6$. There are $6$ such representatives as given below:
\begin{align*}
&x_1=0,\\
& x_1+x_2=1,\\
& x_1+x_2+x_3=1,\\
& x_1+x_2+x_3+x_4=2,\\
& x_1+x_2+x_3+x_4+x_5=2,\\
& x_1+x_2+x_3+x_4+x_5+x_6=3.
\end{align*}

Clearly, two spanned hyperplanes of $Q_n$ in the same
equivalence class contain the same number
of vertices of $Q_n$. So we may say that
an equivalence class of spanned hyperplanes
of $Q_n$ contains $k$ vertices of $Q_n$ if
every hyperplane in this class contains $k$ vertices of $Q_n$.

To state the main result of this section, we need to define the
equivalence classes of $0/1$-polytopes contained in a set of points¡¡ in $\mathbb{R}^n$.
Given a set $\mathcal{S}\subset \mathbb{R}^n$, consider the set of $0/1$-polytopes
of $Q_n$ that are contained in $\mathcal{S}$, denoted by $\mathcal{S}(Q_n)$. Restricting the $0/1$-equivalence relation to the set $\mathcal{S}(Q_n)$
indicates a equivalence relation on $\mathcal{S}(Q_n)$.
More precisely, two $0/1$-polytopes in $\mathcal{S}(Q_n)$
are equivalent if  one can be transformed to the other by a symmetry of $Q_n$.
We call equivalence classes of $0/1$-polytopes  in  $\mathcal{S}(Q_n)$   partial
$0/1$-equivalence classes of $\mathcal{S}$ for the reason that any partial equivalence class of $\mathcal{S}$ is a subset of  a (unique) $0/1$-equivalence class of $Q_n$. Notice that for a $0/1$-polytope $P$ contained in $\mathcal{S}$
and a symmetry $w$, $w(P)$ is not in the partial $0/1$-equivalence class of $P$ when $w(P)$ is not  in $\mathcal{S}(Q_n)$.
Denote by  $\mathcal{P}(\mathcal{S},k)$ the
set of partial $0/1$-equivalence classes of $\mathcal{S}$ with $k$ vertices.
Let  $N_\mathcal{S}(k)$ be the cardinality of $\mathcal{P}(\mathcal{S},k)$.

Let $h(n,k)$ denote the
number of equivalence classes of spanned hyperplanes
of $Q_n$ that contain at least $k$ vertices of $Q_n$.
Assume that  $H_1,H_1,\ldots,H_{h(n,k)}$ are the representatives
of equivalence classes of spanned hyperplanes of
$Q_n$ containing at least $k$ vertices of $Q_n$.
Recall that $\mathcal{H}_n(k)$ denotes the set of $0/1$-equivalence
classes of $Q_n$ with $k$ vertices
that are not full-dimensional. We shall define a map, denoted by $\Phi$, from the (disjoint) union of $\mathcal{P}(H_i,k)$ for $1\leq i\leq h(n,k)$ to $\mathcal{H}_n(k)$. Given a partial $0/1$-equivalence class $\mathcal{P}_i\in \mathcal{P}(H_i,k)$ ($1\leq i\leq h(n,k)$), then we define $\Phi(\mathcal{P}_i)$ to be the (unique) $0/1$-equivalence class in $\mathcal{H}_n(k)$ containing $\mathcal{P}_i$. Then we have the following theorem.

\begin{theorem}\label{cg16}
If $2^{n-2}<k\leq 2^{n-1}$, then the map $\Phi$ is a bijection.
\end{theorem}

\pf
We proceed to  show that $\Phi$ is injective. To this end, we shall prove that for any
two distinct partial $0/1$-equivalence classes $\mathcal{P}_1$
and  $\mathcal{P}_2$ with $k$ vertices, their images, denoted by
$\mathcal{C}_1$ and
$\mathcal{C}_2$, are distinct $0/1$-equivalence classes.
 Assume that $\mathcal{P}_1\in \mathcal{P}(H_i,k)$ and $\mathcal{P}_2\in \mathcal{P}(H_j,k)$, where $1\leq i,j\leq h(n,k)$. Let $P_1$ (resp.  $P_2$) be any $0/1$-polytope in $\mathcal{P}_1$ (resp. $\mathcal{P}_2$). Evidently, $P_1$  (resp.  $P_2$) is a $0/1$-polytope in $\mathcal{C}_1$ (resp. $\mathcal{C}_2$). To prove that
 $\mathcal{C}_1\neq \mathcal{C}_2$, it suffices to show that
   $P_1$ and $P_2$ are not in the same $0/1$-equivalence class.
We have two cases.

Case 1. $i=j$. In this case, it is clear that $P_1$ and $P_2$ are not equivalent.

Case 2. $i\neq j$. Suppose to the contrary  that $P_1$ and $P_2$ are
in the same $0/1$-equivalence class.
Then there exists a symmetry $w\in B_n$ such that $w(P_1)=P_2$.
 Since $2^{n-2}<k\leq 2^{n-1}$, by Theorem \ref{cg14} we see that
$P_1$ and $P_2$ are of dimension $n-1$. Since $P_1$ is contained in $H_i$,  $H_i$
coincides with  the affine space spanned by $P_1$. Similarly, $H_j$ is the affine space spanned by $P_2$. This implies that $w(H_i)=H_j$,
contradicting  the assumption that $H_i$ and $H_j$ belong to distinct equivalence classes of spanned hyperplanes of $Q_n$.
Consequently,  $P_1$ and $P_2$ are not in the same $0/1$-equivalence class.

It remains to show that $\Phi$ is surjective.
For any $\mathcal{C}\in \mathcal{H}_n(k)$, we aim to find a partial
$0/1$-equivalence class such that its image is $\mathcal{C}$. Let $P$ be any $0/1$-polytope in
$\mathcal{C}$. Since $P$ is not full-dimensional, we can find a spanned hyperplane $H$ of $Q_n$ such that $P$ is contained in $H$. It follows that
 $H$ contains at leat $k$ vertices of $Q_n$. Thus there exists a representative $H_j$ ($1\leq j\leq h(n,k)$) such that $H$ is  in the equivalence
  class of $H_j$. Assume that $w(H)=H_j$ for some $w\in B_n$. Then $w(P)$
is contained in $H_j$. It is easily seen that
 under the  map $\Phi$, $\mathcal{C}_i$ is the image of
 the partial $0/1$-equivalence class of $H_j$ containing $w(P)$.
 Thus we conclude that the above map is a bijection. This completes the proof.
\qed

It should also be noted that in the above
proof of Theorem \ref{cg16},  the condition
$2^{n-2}<k\leq 2^{n-1}$  is  required only in Case 2.
When $k< 2^{n-2}$,  the map $\Phi$ may be no longer  an injection.
For the  case $2^{n-3}<k\leq 2^{n-2}$, we will
consider the computation of $H_n(k)$ in Section 8.

As a direct consequence of Theorem \ref{cg16}, we obtain that for $2^{n-2}< k\leq  2^{n-1}$,
\begin{equation}\label{cg15}
H_n(k)=\sum_{i=1}^{h(n,k)}N_{H_i}(k).
\end{equation}
Thus, for $2^{n-2}<k\leq 2^{n-1}$
the computation of $H_n(k)$
is reduced to the determination of the number $N_H(k)$ of
partial $0/1$-equivalence classes of $H$ with $k$ vertices. In the  rest
 of
this section, we shall explain how to compute $N_H(k)$.

For $2^{n-2}<k\leq 2^{n-1}$, let $H$ be a spanned
hyperplane of $Q_n$ containing at least $k$ vertices.  Let $P$ and $P'$ be two distinct 0/1-polytoeps
of $Q_n$ with $k$ vertices that are contained in $H$. Assume that $P$ and $P'$
belong to the same partial  $0/1$-equivalence class of $H$.  Then
there exists a symmetry $w\in B_n$ such that
$w(P)=P'$. It is clear from  Theorem \ref{cg14} that
both $P$ and $P'$ have  dimension
$n-1$. Then $H$ is the affine space spanned by $P$ or $P'$. So we deduce that
$w(H)=H$. Let \[F(H)=\{w\in B_n\colon w(H)=H\}\]
be the stabilizer subgroup of $H$, namely, the subgroup of $B_n$ that  fixes $H$.
So
we have shown that $P$ and $P'$
belong to the same
partial $0/1$-equivalence class of $H$
if and only if one can be transformed to the other by a symmetry in
$F(H)$.

The above fact allows us to use P\'olya's
theorem to compute the number $N_H(k)$ for $2^{n-2}<k\leq 2^{n-1}$.
Denote by $V_n(H)$  the set of vertices of $Q_n$ that are
contained in $H$.  Let us consider the action of $F(H)$ on $V_n(H)$.
 Assume that
 each  vertex in $V_n(H)$ is assigned one of the two colors, say, black and white.
For such a 2-coloring of the vertices in $V_n(H)$,
 consider the black vertices as vertices of a $0/1$-polytope contained in $H$.
  Clearly, for $2^{n-2}< k\leq 2^{n-1}$, this establishes a one-to-one correspondence between partial $0/1$-equivalence classes of $H$ with $k$ vertices and equivalence classes of 2-colorings of the vertices in $V_n(H)$ with $k$ black vertices.

Write $Z_H(z)$  for the cycle index of
$F(H)$, and let  $C_H(u_1,u_2)$ denote the polynomial
obtained from
$Z_H(z)$ by substituting $z_i$ with $u_1^i+u_2^i$.

\begin{theorem}\label{cg111}
 Assume that $2^{n-2}<k\leq 2^{n-1}$, and let $H$ be a spanned
hyperplane of $Q_n$ containing  at least $k$ vertices
of $Q_n$. Then we have
\[N_H(k)=\left[u_1^ku_2^{|V_n(H)|-k}\right]C_H(u_1,u_2).\]
\end{theorem}

We will compute the cycle index $Z_H(z)$ in Section 5 and Section 6.
Section 5 is devoted to the characterization of the stabilizer $F(H)$. In
Section 6, we will give an explicit expression for $Z_H(z)$.

\section{The structure of the stabilizer $F(H)$}

In this section, we aim to
characterize the stabilizer  $F(H)$ for a given spanned
hyperplane $H$  of $Q_n$.

Let \[H\colon a_1x_1+a_2x_2+\cdots+a_nx_n=b\]
be a spanned hyperplane of $Q_n$.
Given $w\in B_n$, let $s(w)$ be the set of entries of $w$ that are assigned the minus sign.
In view of (\ref{cg19}), it is easy to see that
$w(H)$ is of the following form
\begin{equation}\label{linear}
\sum_{i\,\notin s(w)}a_{\pi(i)}x_i+
\sum_{j\,\in s(w)}a_{\pi(j)}(1-x_j)=b,\end{equation}
where $\pi$ is the underlying
permutation of $w$.
The hyperplane $w(H)$ in (\ref{linear}) can be rewritten as
\begin{align}\label{1}
 s(w,1)\cdot a_{\pi(1)}x_1+s(w,2)\cdot a_{\pi(2)}
 x_2+\cdots+s(w,n)\cdot
a_{\pi(n)}x_n=b-\sum_{j\in s(w)}a_{\pi(j)},
\end{align} where
$s(w,j)=-1$ if $j\in s(w)$ and
$s(w,j)=1$ otherwise.

As an example, let \[H\colon x_1-x_2-x_3+2x_4=1\]
be a spanned hyperplane of $Q_4$. Upon the action of the symmetry
$w=(1)(\bar{2}\bar{3})(4)\in
B_4$,  $H$ is transformed into the following hyperplane
 \[x_1+x_2+x_3+2x_4=3.\]

As mentioned in Section 4,
 for
every equivalence class of  spanned hyperplanes of $Q_n$,
we can choose a  representative of the following form
\begin{equation}\label{cg17}
H\colon\, a_1x_1+a_2x_2+\cdots+a_tx_t=b,\end{equation}
where  $a_1\leq a_2\leq
\cdots\leq a_t$ ($t\leq n$), and the coefficients
$a_i$'s and $b$  are positive integers.
Note that this observation also follows from (\ref{1}).
From now on, we shall restrict our attention only to spanned hyperplanes
of $Q_n$ of the form as in (\ref{cg17}).
 The following definition is required for the determination of $F(H)$.

\begin{definition}\label{defi}
Let $H$ be a spanned hyperplane of the form as in (\ref{cg17}). The   type of $H$ is defined to be a vector $\alpha=(\alpha_1,\alpha_2,\ldots,\alpha_\ell)$,
where $\alpha_i$ is the multiplicity of $i$ occurring in the  set
$\{a_1,a_2,\ldots,a_t\}$.
\end{definition}

For example, let
\begin{equation}\label{hyperplane}
H\colon x_1+x_2+2x_3+2x_4+3x_5=4\end{equation}
be a spanned hyperplane of $Q_5$. Then, the type of $H$ is
$\alpha=(\alpha_1,\alpha_2,\alpha_3)=(2,2,1).$

For positive integers $i$ and $j$ such that $i\leq j$, let
$[i,j]$ denote the interval $\{i,i+1,\ldots,j\}$.
Let
$\alpha=(\alpha_1,\alpha_2,\ldots,\alpha_\ell)$ be the type of
a spanned hyperplane. Under the
assumption that $\alpha_0=0$,  the following set
\begin{equation}\label{partition}
\left\{\left[\alpha_1+\cdots+\alpha_{i-1}+1,\alpha_1+\cdots+
\alpha_{i-1}+\alpha_i\right]
\colon 1\leq i\leq \ell\right\}\end{equation}  is a partition
of the set $\{1,2,\ldots,t\}$. For example,
 let $\alpha=(2,2,1)$. Then the corresponding partition  is $\{\{1,2\},\{3,4\},\{5\}\}$.

Since   (\ref{partition}) is a partition of $\{1,2,\ldots,t\}$,
we can define the corresponding Young subgroup $S_\alpha$
of the permutation group on $\{1,2,\ldots,t\}$, namely,
\begin{equation}\label{twin1}
S_\alpha=S_{\alpha_1}\times S_{\alpha_2}\times\cdots\times S_{\alpha_\ell},\end{equation}
where $\times$ denotes the direct
product of groups, and for $i=1,2,\ldots,\ell$, $S_{\alpha_i}$ is the permutation group on the interval
\begin{equation}\label{set}
\left[\alpha_1+\cdots+\alpha_{i-1}+1,\alpha_1+\cdots+
\alpha_{i-1}+\alpha_i\right].\end{equation}
Let
\begin{equation}\label{twin2}
\overline{S}_\alpha=\overline{S}_{\alpha_1}\times
\overline{S}_{\alpha_2}\times\cdots\times
\overline{S}_{\alpha_\ell},\end{equation}
where
$\overline{S}_{\alpha_i}$ is the set of signed permutations on
the interval (\ref{set}) with
all elements assigned the minus sign.
Define
\begin{equation}\label{Fixgroup}
P(H)=\left\{\begin{array}{lll}
S_\alpha, & \ \  \mathrm{if} \ \ \sum_{i=1}^t a_i\neq 2b,\\&\\
S_\alpha\bigcup \overline{S}_\alpha, & \ \ \mathrm{if}\ \
\sum_{i=1}^t a_i=2b.
\end{array}\right.
\end{equation}
The following theorem  gives a characterization of
the stabilizer of a spanned hyperplane.

\begin{theorem}\label{cg4}
Let $H\colon a_1x_1+a_2x_2+\cdots +a_tx_t=b$ be
a spanned hyperplane of $Q_n$. Then
\[F(H)=P(H)\times B_{n,t},\] where $B_{n,t}$ is the group
of all signed permutations on the interval $[t+1,n]$.
\end{theorem}

\pf
Assume that $w\in F(H)$ and $\pi$ is the underlying permutation of $w$. We  aim to show that $w\in P(H)\times B_{n,t}$.
Consider the expression of $w(H)$ as in (\ref{1}), that is,
\begin{equation}\label{cg101}
 s(w,1)\cdot a_{\pi(1)}x_1+s(w,2)\cdot a_{\pi(2)}
 x_2+\cdots+s(w,n)\cdot
a_{\pi(n)}x_n=b-\sum_{j\in s(w)}a_{\pi(j)}.
\end{equation}
We claim  that $s(w,j)$  are either all positive or all negative for $1\leq j\leq t$.  Suppose otherwise that there exist $1\leq i,j\leq t$ ($i\neq j$) such  that $s(w,i)>0$ and $s(w,j)<0$. Since the $a_i$'s are all positive,  we see that the coefficients $s(w,i)a_{\pi(i)}$ and $s(w,j)a_{\pi(j)}$ for the hyperplane $w(H)$ have opposite signs. This  implies that $w(H)$ and $H$ are distinct, which contradicts the assumption that $w$ fixes $H$.   We now have the following two cases.

Case 1.  The signs $s(w,j)$ are all positive for $1\leq j\leq t$.  In this case, since $w(H)=H$
it is clear that $w(H)$ is of the following form
\[a_{\pi(1)}x_1+a_{\pi(2)}x_2+\cdots+a_{\pi(t)}x_t=b,\] where $a_{\pi(j)}=a_j$ for $1\leq j\leq t$. Hence we deduce that, for any $1\leq j\leq t$, $\pi(j)$ is in the interval $\left[\alpha_1+\cdots+\alpha_{i-1}+1,\alpha_1+\cdots+
\alpha_{i-1}+\alpha_i\right]$ that  contains the element $j$.
Thus we obtain that $w\in S_\alpha\times B_{n,t}$.

Case 2. The signs  $s(w,j)$ are all negative for $1\leq j\leq t$.  In this case,
 we see that $w(H)$ is of the following form
\[-a_{\pi(1)}x_1-a_{\pi(2)}x_2-\cdots-a_{\pi(t)}x_t=b-(a_1+\cdots+a_t).\]
Since $w(H)=H$,
 we have $a_{\pi(j)}=a_j$ for $1\leq j\leq t$  and $b-(a_1+\cdots+a_t)=-b$.
Thus we obtain  $w\in \overline{S}_\alpha\times B_{n,t}$.
Combining the above two cases, we conclude that $w\in P(H)\times B_{n,t}$.

On the other hand, from the expression (\ref{cg101})
for $w(H)$, it is not difficult to check that  every symmetry $w$ in $P(H)\times B_{n,t}$ fixes $H$. This completes the proof. \qed

As will been seen in Section 6, for the purpose of computing the
cycle index $Z_H(z)$ with respect to a spanned hyperplane $H$ of $Q_n$,
it is often necessary to consider the structure of the subgroup $P(H)$ of $F(H)$.
We sometimes write a symmetry $\pi\in P(H)$ as a product form $\pi=\pi_1\pi_2\cdots\pi_\ell$, which means that for $i=1,2\ldots,\ell$, $\pi_i\in S_{\alpha_i}$ if $\pi\in S_\alpha$, and  $\pi_i\in \overline{S}_{\alpha_i}$ if $\pi\in \overline{S}_\alpha$, where $\alpha$
is the type of $H$.
We conclude this section with the following proposition, which will be required for the computation of $Z_H(z)$ in Section 6.

\begin{prop}\label{prop1}
Let $H$ be a spanned hyperplane of $Q_n$ of type $\alpha$. Let $\pi=\pi_1\pi_2\cdots\pi_\ell$ and $\pi'=\pi_1'\pi_2'\cdots\pi_\ell'$ be two symmetries in  $P(H)$, and assume that both $\pi$ and $\pi'$ are either in
$ S_\alpha$ or in $\overline{S}_\alpha$. If $\pi_i$ and $\pi_i'$ have the same cycle type for $1\leq i\leq \ell$, then $\pi$ and $\pi'$ are in the same conjugacy class of $P(H)$.
\end{prop}

\pf To prove that $\pi$ and $\pi'$ are conjugate in $P(H)$, it suffices
to show that there exists a symmetry $w\in P(H)$ such that $\pi=w\pi'w^{-1}$. First, we consider the case when both $\pi$ and $\pi'$ are in $S_\alpha$. Since $\pi_i$ and $\pi_i'$ are
of the same cycle type,  they are in the same
conjugacy class. So there is a permutation  $w_i\in S_{\alpha_i}$ such that $\pi_i=w_i\pi_i'w_i^{-1}$. It follows that $\pi=(w_1\pi_1'w_1^{-1})\cdots(w_\ell\pi_\ell'w_\ell^{-1})=w\pi'w^{-1}$, where $w=w_1\cdots w_\ell\in S_\alpha$. This implies that $\pi$ and $\pi'$ are conjugate in $P(H)$.

It remains to consider the case when both $\pi$ and $\pi'$ are in
$\overline{S}_\alpha$.
Let $\pi_0$ (resp. $\pi'_0$) be the underlying permutation of $\pi$ (resp. $\pi'$).
Then there is a symmetry $w\in S_\alpha$ such that
$\pi_0=w\pi'_0w^{-1}$. We claim that $\pi=w\pi'w^{-1}$.  Indeed,
 it is enough to show that $\pi(x_1,x_2,\ldots,x_t)
 =w\pi'w^{-1}(x_1,x_2,\ldots,x_t)$ for  any point $(x_1,x_2,\ldots,x_t)$ in $\mathbb{R}^t$.
Assume that  $\pi(x_1,x_2,\ldots,x_t)=(y_1,y_2,\ldots,y_t)$ and $w\pi'w^{-1}(x_1,x_2,\ldots,x_t)=(z_1,z_2,\ldots,z_t)$.
Since all elements of $\pi$ are assigned the minus sign, we obtain from (\ref{cg19})  that $y_i=1-x_{\pi_0(i)}$ for $1\leq i\leq t$. On the other hand, using (\ref{cg19}), it is not hard to check that $z_i=1-x_{w^{-1}\pi'_0w(i)}$ for $1\leq i\leq t$.
Since $\pi_0=w\pi'_0w^{-1}$, we deduce that $\pi_0(i)=w^{-1}\pi'_0w(i)$.
Therefore, we have $y_i=z_i$ for  $1\leq i\leq t$. So the claim is justified.
This completes the proof.
\qed

\section{The computation of $Z_H(z)$}

 In this section, we shall
derive a formula  for the cycle index $Z_H(z)$ for a spanned hyperplane $H$ of $Q_n$. It turns out that $Z_H(z)$ depends only on the cycle structures of the symmetries in the subgroup $P(H)$ of $F(H)$.

Let
\begin{equation}\label{dfod}
H\colon\, a_1x_1+a_2x_2+\cdots+a_tx_t=b\end{equation}
be a spanned hyperplane of $Q_n$.
Recall that $V_n(H)$ is the set of vertices of $Q_n$ contained in $H$. To compute the cycle index $Z_H(z)$, we need to determine the cycle structures of permutations on $V_n(H)$ induced by the symmetries in $F(H)$.
By Theorem \ref{cg4}, each symmetry in $F(H)$ can be
written uniquely as a product $\pi w$,
where $\pi\in P(H)$ and $w\in
B_{n,t}$. We shall define two group actions for the subgroups $P(H)$ and $B_{n,t}$,  and shall derive an expression for the cycle type of the permutation on $V_n(H)$ induced by $\pi w$ in terms of the cycle types  of the permutations induced by $\pi$ and $w$.

Let $H$ be a spanned hyperplane of $Q_n$ as in (\ref{dfod}). However, to define the action for $P(H)$,  we shall consider $H$ as a hyperplane in $\mathbb{R}^t$.  Denoted by $V_t(H)$ the set of vertices of $Q_t$ that are contained in $H$, namely, \[V_t(H)=\{(x_1,x_2,\ldots,x_t)\in V_t\colon a_1x_1+a_2x_2+\cdots+a_tx_t=b\}.\]
Since the vertices of $Q_n$ contained in $H$ span a
hyperplane in $\mathbb{R}^n$,  it can be seen that
the vertices in $V_t(H)$ span a hyperplane in $\mathbb{R}^t$.
 Since $H$ is considered as a hyperplane in $\mathbb{R}^t$,
  we deduce that $H$ is a spanned hyperplane of $Q_t$.
Setting $n=t$ in Theorem \ref{cg4}, it follows that the stabilizer
of $H$ is  $P(H)$. Therefore, $P(H)$ stabilizes the set $V_t(H)$.
So any symmetry in $P(H)$ induces a permutation on  $V_t(H)$.

We also need  an action of the group
$B_{n,t}$ on the set of vertices of $Q_{n-t}$.  Assume that  $w\in B_{n,t}$, namely,
  $w$ is a signed permutation on the interval $[t+1,n]$.
Subtracting each element of $w$ by $t$,
we get a signed permutation on $[1,n-t]$.  In this way, each signed permutation in $B_{n,t}$  corresponds to a symmetry of  $Q_{n-t}$. Hence $B_{n,t}$ is isomorphic to the group $B_{n-t}$  of symmetries of $Q_{n-t}$. This
 leads to an action of the group $B_{n,t}$ on $V_{n-t}$.

Let $\pi\in P(H)$ and  $w\in B_{n,t}$. Recall that, for an element $g$ in a group $G$ acting on a finite set $X$, $c(g)$ denotes the cycle type of the permutation on $X$ induced by $g$, which is written as a multiset $\{1^{c_1}, 2^{c_2}, \ldots\}$.
In this notation, $c(\pi)$ (resp. $c(w)$) represents the  cycle type of the permutation  on $V_t(H)$ (resp. $V_{n-t}$) induced by $\pi$ (resp. $w$).
 The following lemma gives an expression for the cycle type $c(w\pi)$ of the induced permutation of $\pi w$ on $V_n(H)$  in terms of
 the cycle types $c(\pi)$ and $c(w)$.

\begin{lemma}\label{lemma3}
Let $H\colon\, a_1x_1+a_2x_2+\cdots+a_tx_t=b$ be a spanned hyperplane of $Q_n$, and $\pi w$ be a symmetry in $F(H)$,
where $\pi\in P(H)$ and $w\in
B_{n,t}$. Assume that $c(\pi)=\{1^{m_1},2^{m_2},\ldots\}$ and $c(w)=\{1^{c_1},2^{c_2},\ldots\}$. Then we have
\begin{equation}\label{cg24}
c(\pi w)=\bigcup_{i\geq 1}\bigcup_{j\geq
1}\left\{\left(\mathrm{lcm}(i,j)\right)^{\frac{ijm_ic_j}
{\mathrm{lcm}(i,j)}}\right\},\end{equation}
where  $\bigcup$ denotes the disjoint union of multisets, and $\mathrm{lcm}(i,j)$ denotes the least common multiple of integers
$i$ and  $j$.
\end{lemma}

\pf
Clearly, each  vertex in  $V_n(H)$ can be expressed as  a vector of the following
form
\[(x_1,\ldots,x_{t},y_1,\ldots,y_{n-t}),\] where $(x_1,\ldots,x_{t})
$ is a vertex in $V_t(H)$ and $(y_1,\ldots,y_{n-t})$ is a vertex of $Q_{n-t}$.  Assume that $|V_t(H)|=n_0$. Let $V_t(H)=\{u_1,u_2,
\ldots,u_{n_0}\}$ and $Q_{n-t}=\{v_1,v_2,\ldots,v_{2^{n-t}}\}$. Then each vertex in $V_n(H)$ can be expressed  as an ordered pair $(u_i,v_j)$,
where $1\leq i\leq n_0$ and $1\leq j\leq 2^{n-t}$.

Let $C_i=({s_1},\ldots,{s_i})$ be
an $i$-cycle of the permutation on $V_t(H)$ induced by
$\pi$, that is, $C_i$ maps the vertex $u_{s_k}$ to the vertex $u_{s_{k+1}}$ if $1\leq k\leq i-1$, and to the vertex $u_{s_1}$ if $k=i$.
Similarly,  let $C_j=({t_1},\ldots,t_j)$ be a $j$-cycle of the
permutation on $V_{n-t}$ induced by $w$, that is, $C_j$ maps the vertex $v_{t_m}$ to the vertex $v_{t_{m+1}}$ if $1\leq m\leq j-1$, and to the vertex $v_{t_1}$ if $m=j$. Define  the direct product of $C_i$ and $C_j$, denoted $C_i\times C_j$,
 to be the
permutation on the subset
$\{(u_{s_k},v_{t_m})\colon 1\leq k\leq i,\ 1\leq
m\leq j\}$ of $V_n(H)$ such that
\[C_i\times
C_j(u_{s_k},v_{t_m})
=(C_i(u_{s_k}),C_j(v_{t_m})).\]
It is not hard to check  that the cycle type of $C_i\times C_j$
is
\[\left\{\left(\mathrm{lcm}(i,j)\right)^
{\frac{ij}{\mathrm{lcm}(i,j)}}\right\}.\]

Note that  the induced  permutation of $\pi w$ on $V_n(H)$
 is the product of $C_i\times C_j$,
where $C_i$ (resp. $C_j$) runs over the cycles of the permutation on $V_t(H)$ (resp. $V_{n-t}$) induced by $\pi$ (resp. $w$). Thus the cycle type of the induced permutation of $\pi w$ on $V_n(H)$
is  given by (\ref{cg24}). This completes the proof.  \qed

Before presenting a formula for the cycle index $Z_H(z)$, we need to introduce some notation. Assume that $\pi$ is a symmetry in $P(H)$ such that
 the cycle type of the induced permutation of $\pi$ is
\[c(\pi)=\{1^{m_1},2^{m_2}, \ldots\}.\] For  $j\geq 1$, we define
\begin{equation}\label{fdefinition}
f_\pi(z_j)=\prod_{i\geq 1}(z_{\mathrm{lcm}(i,j)})^{\frac{ijm_i}{\mathrm{lcm}(i,j)}}.\end{equation}
Let
\begin{equation}\label{fdefinition2}
f_\pi(z)=(f_\pi(z_1),f_\pi(z_2),\ldots).\end{equation}
We have the following proposition.

\begin{prop}\label{prop10}
Let $H$ be a spanned hyperplane of $Q_n$ with type $\alpha$. Let $\pi=\pi_1\pi_2\cdots\pi_\ell$ and $\pi'=\pi_1'\pi_2'\cdots\pi_\ell'$ be two symmetries in  $P(H)$. Assume that both $\pi$ and $\pi'$ are either in
$ S_\alpha$ or in $\overline{S}_\alpha$. If $\pi_i$ and $\pi_i'$ have the same cycle type for $1\leq i\leq \ell$, then  $f_{\pi}(z)=f_{\pi'}(z)$.
\end{prop}

\pf It follows from Proposition \ref{prop1} that $\pi$ and $\pi'$ are conjugate in $P(H)$.
Since $P(H)$ acts on $V_t(H)$, the permutations on $V_t(H)$ induced by $\pi$ and $\pi'$ are conjugate. So they have the same cycle type, i.e.,
$c(\pi)=c(\pi')$.  Since $f_{\pi}(z)$ depends only on $c(\pi)$, we see that $f_{\pi}(z)=f_{\pi'}(z)$.
This completes the proof. \qed

We now give an overview of  some notation related to integer partitions. We shall write a partition $\lambda$ of a positive integer $n$, denoted by $\lambda\vdash n$,   in the multiset form, that is, write $\lambda=\{1^{m_1},2^{m_2},\ldots\}$, where $m_i$ is the number of parts of $\lambda$ of size $i$. Denote by $\ell(\lambda)$  the number of parts of $\lambda$, that is, $\ell(\lambda)=m_1+m_2+\cdots$. For a partition $\lambda=\{1^{m_1},2^{m_2},\ldots\}$,
let \[z_\lambda=1^{m_1}
m_1! 2^{m_2}m_2!\cdots.\]
For two partitions $\lambda$ and $\mu$, define $\lambda \cup \mu$
to be the partition obtained by joining the parts of $\lambda$ and
$\mu$ together. For example, for $\lambda=\{1,2\}$ and $\mu=\{1^2,3\}$,
then $\lambda \cup \mu=\{1^3,2,3\}$.

Let $H\colon a_1x_1+a_2x_2+\cdots+a_tx_t=b$ be a spanned hyperplane of $Q_n$, whose  type is
$\alpha=(\alpha_1,\alpha_2,\ldots,\alpha_\ell)$. Assume that  $\mu=\mu^1\cup\cdots\cup\mu^\ell$ is a partition of $t$,
where $\mu^i\vdash \alpha_i$ for $1\leq i\leq \ell$. We can write  $f_\mu(z)$ (resp. $\overline{f}_\mu(z)$)  for $f_{\pi}(z)$, where $\pi=\pi_1\pi_2\cdots\pi_\ell$ is any symmetry in $S_\alpha$
(resp. $\overline{S}_\alpha$) such that $\pi_i$  has cycle type $\mu^i$ for $1\leq i\leq \ell$. By Proposition \ref{prop10}, the
functions $f_\mu(z)$ and $\overline{f}_\mu(z)$ are well defined.
We can now give a formula for the cycle index $Z_H(z)$.

\begin{theorem}\label{cg6}
Let $H\colon a_1x_1+a_2x_2+\cdots+a_tx_t=b$ be a spanned hyperplane of $Q_n$. Assume that $H$ has type
$\alpha=(\alpha_1,\alpha_2,\ldots,\alpha_\ell)$. Then we have
\begin{equation}\label{mainone}
Z_H(z)=\frac{1}{2^{\delta(H)}}\sum_{(\mu^1,\ldots,\mu^\ell)}
\prod_{i=1}^\ell z_{\mu^i}^{-1}\left(Z_{n-t}(f_\mu(z))+
\delta(H)Z_{n-t}(\overline{f}_\mu(z))\right),\end{equation}
where $\mu^i\vdash \alpha_i$,
$\mu=\mu^1\cup\cdots\cup\mu^\ell$,
$\delta(H)=1$ if $\sum_{i=1}^t a_i= 2b$ and
$\delta(H)=0$ otherwise.
\end{theorem}

\pf
Let  $\pi\in P(H)$ and $w\in B_{n,t}$. Assume that $c(w)=\{1^{c_1},2^{c_2},\ldots\}$. By Lemma \ref{lemma3}, we have
\begin{equation}\label{cg112}
{z}^{c(\pi\cdot
w)}=f_\pi(z_1)^{c_1}f_\pi(z_2)^{c_2}\cdots.
\end{equation}
From (\ref{cg22}) and (\ref{cg112}),  we deduce that
\begin{align*} \sum_{\pi w}{z}^{c(\pi\cdot
w)}&=\sum_{w}f_\pi(z_1)^{c_1}f_\pi(z_2)^{c_2}\cdots\\
&=(n-t)!2^{n-t}Z_{n-t}(f_\pi(z_1),f_\pi(z_2),\ldots)\\[5pt]
&=(n-t)!2^{n-t}Z_{n-t}(f_\pi({z})),
\end{align*}
where $w$ runs over the signed permutations in $B_{n,t}$. Thus
\begin{equation}\label{cg25}
\begin{split}
Z_H({z})&=\frac{1}{|F(H)|}\sum_{\pi w\in
F(H)}{z}^{c(\pi w)}\\[5pt]
&=\frac{1}{|F(H)|}\sum_{\pi\in
P(H)}(n-t)!2^{n-t}Z_{n-t}(f_\pi({z}))\\[5pt]
&=\frac{(n-t)!2^{n-t}}{|F(H)|}\Bigg(\sum_{\pi\in
S_\alpha}Z_{n-t}(f_\pi({z}))+
\delta(H)\sum_{\pi'\in
\overline{S}_\alpha}Z_{n-t}(f_{\pi'}({z}))\Bigg),
\end{split}
\end{equation}
where $\delta(H)=1$ if $\sum_{i=1}^t a_i= 2b$ and
$\delta(H)=0$ otherwise.

Recall that for any given partition $\nu\vdash m$,
there are $\frac{m!}{z_\nu}$ permutations on $\{1,2,\ldots,m\}$ such that their  cycle type is
$\nu$, see Stanley \cite[Proposition 1.3.2]{Stanley}.
So the number of symmetries $\pi=\pi_1\pi_2\cdots\pi_\ell$ in $S_\alpha$ (or, $\overline{S}_\alpha$)  such that for $i=1,2\ldots,\ell$, $\pi_i$ has cycle type $\mu^i$ is equal to
\begin{equation}\label{product}
\prod_{i=1}^\ell\frac{
\alpha_i!}{z_{\mu^i}}.
\end{equation}
Combining (\ref{cg25}),  (\ref{product}) and Proposition \ref{prop10}, we obtain that
\begin{equation}\label{cg26}
Z_H({z})
=\frac{(n-t)!2^{n-t}}{|F(H)|}\sum_{(\mu^1,\ldots,\mu^\ell)}
\prod_{i=1}^\ell\frac{
\alpha_i!}{z_{\mu^i}}\big(Z_{n-t}(f_\mu({z}))+
\delta(H)Z_{n-t}(\overline{f}_\mu({z}))\big),
\end{equation}
where $\mu^i\vdash \alpha_i$, and
$\mu=\mu^1\cup\cdots\cup\mu^\ell$.

Since
\begin{equation}\label{cg27}
|F(H)|=(n-t)!2^{n-t+\delta(H)}\prod_{i=1}^\ell
\alpha_i!,\end{equation}
by substituting (\ref{cg27}) into  (\ref{cg26}), we are led to  (\ref{mainone}). This completes the proof.
\qed

By Theorem \ref{cg6}, the cycle index $Z_H({z})$
depends only on $f_{\pi}(z)$ for $\pi\in P(H)$.
In view of (\ref{fdefinition}), we see that $f_{\pi}(z)$ depends only
 on  $c(\pi)$.
Assume that
$c(\pi)=\{1^{m_1},2^{m_2},\ldots\}$. By Theorem \ref{cg7}, we  have
\begin{equation}\label{mobius}
m_i=\frac{1}{i}\sum_{j|i}\mu(i/j)\psi(\pi^j),\end{equation}
where  $\psi(\pi^j)$ is the number of vertices in $V_t(H)$ that are fixed by  $\pi^j$.
The following theorem gives a formula for $\psi(\pi)$, from which
 $\psi(\pi^j)$ is easily determined.

\begin{theorem}\label{30}
Let $H\colon a_1x_1+a_2x_2+\cdots+a_tx_t=b$ be a spanned hyperplane of $Q_n$.
Assume that  $\pi=\pi_1\pi_2\cdots \pi_\ell$ is a symmetry in $P(H)$ such that $\pi_i$ has cycle type  $\mu^i=\{1^{m_{i1}},2^{m_{i2}},\ldots\}$ for $i=1,2,\ldots,\ell$. Then
\begin{equation}\label{2}
\psi(\pi)=\left\{\begin{array}{lll}
\left[x^b\right]\prod_{i=1}^\ell\prod_{j\geq 1}(1+x^{ij})^{m_{ij}}, & \ \  \mathrm{if} \ \ \pi\in S_\alpha,\\&\\
\chi(\mu)2^{\ell(\mu)}, & \ \ \mathrm{if}\ \  \pi\in \overline{S}_\alpha,
\end{array}\right.
\end{equation}
where $\mu=\mu^1\cup\cdots\cup\mu^\ell$, $\chi(\mu)=1$ if $\mu$
has no odd parts and $\chi(\mu)=0$ otherwise.
\end{theorem}

Before we present the proof of the above theorem, we need to define $0/1$-labelings of a symmetry $\pi\in P(H)$ for the purpose of characterizing the vertices of $Q_t$  fixed by $\pi$.  Let $\pi$ be a symmetry in $P(H)$. A $0/1$-labeling of $\pi$ is a labeling of the cycles of $\pi$
such that each cycle of $\pi$ is assigned one of the two numbers $0$ and $1$.

\noindent\textit{Proof of Theorem \ref{30}.}
We first consider the case  when  $\pi$ is in $S_\alpha$.
It is easy to observe that,
a vertex  $v=(v_1,v_2,\ldots,v_t)$ of $Q_t$ is a fixed point of $\pi$,
that is,  $\pi(v)=v$
if and only if, for each $i$-cycle $({j_1},{j_2},\ldots,{j_i})$ of $\pi$ and
for any entry of $v$ corresponding to $({j_1},{j_2},\ldots,{j_i})$, we have
\[v_{{j_1}}=v_{{j_2}}=\cdots=v_{{j_i}}\] (or, more precisely,
$v_{{j_1}}=v_{{j_2}}=\cdots=v_{{j_i}}=0$ or
$v_{{j_1}}=v_{{j_2}}=\cdots=v_{{j_i}}=1$).
The above characterization enables us to establish a one-to-one
correspondence between $0/1$-labelings of $\pi$ and the vertices
of $Q_t$ fixed by $\pi$, that is, for any given  $0/1$-labeling of $\pi$, we can define a  vertex $v=(v_1,v_2,\ldots,v_t)$ of $Q_t$ fixed by $\pi$ such that $v_i=0$ ($1\leq i\leq t$) if and only if the
cycle of $\pi$ containing $i$ is assigned  $0$.
Moreover, if the vertex $v=(v_1,v_2,\ldots,v_t)$ corresponding  to a $0/1$-labeling of $\pi$ is in $V_t(H)$, that is, $a_1v_1+a_2v_2+\cdots+a_tv_t=b$, then we have
\begin{equation}\label{condition1}
b_1+2b_2+\cdots
+\ell b_\ell=b,\end{equation} where $b_i$
 ($1\leq i\leq \ell$) is the sum of the lengths of
cycles of $\pi_i$ which are labeled  $1$.
It can be easily  deduced  that the number of $0/1$-labelings of $\pi$ satisfying (\ref{condition1}) is
\[\psi(\pi)=\left[x^b\right]\prod_{i=1}^\ell\prod_{j\geq 1}\left(1+x^{ij}\right)^{m_{ij}}.\]

We now consider the case when $\pi$ is in $\overline{S}_\alpha$.
As in the previous case, it can be seen that
 a vertex $v=(v_1,v_2,\ldots,v_t)$ of $Q_t$ is fixed by $\pi$
if and only if,  for any (signed)
$i$-cycle
$(\overline{j_1},\overline{j_2},\ldots,\overline{j_i})$
of $\pi$, the following relation holds
\begin{equation}\label{condition2}
(v_{{j_1}},v_{{j_2}},\ldots,v_{{j_i}})
=(1-v_{{j_2}},1-v_{{j_3}},\ldots,1-v_{{j_1}}).\end{equation}
Consequently, if a vertex $v=(v_1,v_2,\ldots,v_t)$
 of $Q_t$ is fixed by $\pi$, then, for any (signed) $i$-cycle
$(\overline{j_1},\overline{j_2},\ldots,\overline{j_i})$
of $\pi$, the vector $(v_{{j_1}},v_{{j_2}},\ldots,v_{{j_i}})$ is
 either $(0,1,\ldots,0,1)$ or $(1,0,\ldots,1,0)$.
This implies that  $\pi$ does not have any fixed point if
 $\pi$  has an odd cycle.

 We now assume that $\pi$ has only even (signed) cycles. In this case, we see that the number of vertices of $Q_t$  fixed by $\pi$ is equal to $2^{\ell(\mu)}$. To prove $\psi(\pi)=2^{\ell(\mu)}$, we need to demonstrate  that any vertex of $Q_t$ fixed by $\pi$ is in $V_t(H)$.
Let  $v=(v_1,v_2,\ldots,v_t)$ be any vertex of $Q_t$  fixed by $\pi$.
Using the fact that for each (signed)  cycle  $(\overline{j_1},\overline{j_2},\ldots,\overline{j_i})$ of $\pi$, the vector $(v_{{j_1}},\ldots,v_{{j_i}})$ is $(1,0,\ldots,1,0)$ or $(0,1,\ldots,0,1)$, and applying the relation
\begin{equation*}
a_1+\cdots+a_t=2b,\end{equation*} we  deduce  that $a_1v_1+a_2v_2+\cdots+a_tv_t=b$. Hence the vertex $v$ is in $V_t(H)$. This completes the proof. \qed

By Theorem \ref{30}, we can compute  $\psi(\pi^j)$.
Let  $\pi=\pi_1\pi_2\cdots \pi_\ell\in P(H)$, where $\pi_i$ has cycle type  $\mu^i=\{1^{m_{i1}},2^{m_{i2}},\ldots\}$. Clearly,  $\pi^j=\pi_1^j\pi_2^j\cdots\pi_\ell^j$. Let
$\mathrm{gcd}(i,j)$ denote the greatest common divisor of $i$ and  $j$.
As is easily checked,  the cycle type of  $\pi_i^j$ ($1\leq i\leq \ell$) is
\[\left\{1^{m_{i1}},
\mathrm{gcd}(2,j)^{\frac{2m_{i2}}{\mathrm{gcd}
(2,j)}}, \mathrm{gcd}(3,j)^{\frac{3m_{i3}}{\mathrm{gcd}
(3,j)}}, \ldots\right\}.\]
 To apply Theorem \ref{30},  it is still necessary to determine whether the symmetry $\pi^j$ belongs to
   $S_\alpha$ or  $\overline{S}_\alpha$.
It can be seen  that if $\pi$ is in $S_\alpha$    or $\pi$ is in $\overline{S}_\alpha$ and $j$ is even, then $\pi^j$ belongs to $S_\alpha$.
Similarly,  if $\pi$ is in $\overline{S}_\alpha$ and $j$ is odd, then $\pi^j$
belongs to $\overline{S}_\alpha$.

\section{$F_n(k)$ for $n=4,5,6$ and $2^{n-2}< k\leq 2^{n-1}$}

This section is devoted to the computation of $F_n(k)$ for $n=4,5,6$ and $2^{n-2}< k\leq 2^{n-1}$. This requires
the cycle indices $Z_H(z)$ for  spanned hyperplanes of $Q_n$ for $n=4,5,6$ that contain more than $2^{n-2}$ vertices of $Q_n$.

Let $H_1,H_2,\ldots,H_{h(n,k)}$ be the representatives of equivalence classes of spanned hyperplanes of $Q_n$ containing at least $k$ vertices. When $2^{n-2}< k\leq 2^{n-1}$, combining relation (\ref{cg10}), Theorem \ref{cg16} and Theorem \ref{cg111}, we deduce that
\begin{equation}\label{cg103}
\begin{split}
F_n(k)&=A_n(k)-H_n(k)\\
&=A_n(k)-\sum_{i=1}^{h(n,k)}N_{H_i}(k)\\
&=A_n(k)-\sum_{i=1}^{h(n,k)}\left[u_1^{k}u_2^{|V_n(H_i)|-k}\right]C_{H_i}(z_1,z_2).
\end{split}
\end{equation}

We start with the  computation of $F_4(k)$ for $k=5, 6, 7, 8$. Observing that  $F_4(k)=0$ for $k<5$, this gives
the enumeration of  full-dimensional $0/1$-equivalence classes of $Q_4$.
For brevity, we use $H_n^t$ ($t\leq n$) to denote the following hyperplane in $\mathbb{R}^n$ \[x_1+x_2+\cdots+x_t=\left\lfloor t/2\right\rfloor.\]
In this notation, representatives  of equivalence classes of spanned hyperplanes of $Q_4$  containing
more than  $4$ vertices of $Q_4$ are as follows
\begin{align*}
&H_4^1\colon x_1=0,\\
&H_4^2\colon x_1+x_2=1,\\
&H_4^3\colon x_1+x_2+x_3=1,\\
&H_4^4\colon x_1+x_2+x_3+x_4=2.
\end{align*}
Employing the techniques in Section 6, we  obtain the cycle indices
$Z_{H_4^1}(z)$ and $Z_{H_4^2}(z)$ as given below:
\begin{align*}
Z_{H_4^1}(z)&=Z_3(z),\\[5pt]
Z_{H_4^2}(z)&= \frac{1}{16}\left(\begin{array}{l}
9z_2^4+4z_4^2+2z_1^4z_2^2+z_1^8\end{array}\right).
\end{align*}
For the remaining two hyperplanes  $H=H_4^3$ and $H_4^4$, it can be  checked that $N_H(k)=1$ for $k=5,6$.
Thus, from (\ref{cg103}) we can determine $F_4(k)$ for $k=5,6,7,8$. These values  are given in Table \ref{t11}, which agree
with the  results computed  by  Aichholzer \cite{Aic2}.
\begin{table}[h,t]
\begin{center}
\begin{tabular}{|l|l|l|l|l|l|l|l|l|l|l|l|l|l|}\hline
$k$&5&6&7&8\\\hline
$H_4^1$&3&3&1&1\\\hline
$H_4^2$&5&5&1&1\\\hline
$H_4^3$&1&1&&\\\hline
$H_4^4$&1&1&&\\\hline
$F_4(k)$&17&40&54&72\\\hline
\end{tabular}
\end{center}
\caption{$F_4(k)$ for $k=5,6,7,8$.}\label{t11}
\end{table}

We now compute $F_5(k)$ for $8< k\leq
16$. Representatives  of equivalence classes of spanned hyperplanes of $Q_5$  containing
more than  $8$ vertices of $Q_5$ are $H_5^1,H_5^2,H_5^3,H_5^4,H_5^5$. By utilizing the  the techniques in Section 6, we  obtain that
\begin{align*}
Z_{H_5^1}({z})&=Z_4({z}),\\[5pt]
Z_{H_5^2}({z})&=\frac{1}{96}\left(\begin{array}{l}z_1^{16}+6z_1^8z_2^4+33z_2^8+8z_1^4z_3^4+24z_4^4+
24z_2^2z_6^2\end{array}\right),\\[5pt]
Z_{H_5^3}({z})&=\frac{1}{48}\left(\begin{array}{l}12z_2^6+8z_4^3+2z_1^6z_2^3+
z_1^{12}+6z_1^2z_2^5+
3z_1^4z_2^4+ 6z_6^2+4z_{12}+4z_3^2z_6+2z_3^4\end{array}\right),\\[5pt]
Z_{H_5^4}({z})&=\frac{1}{96}\left(\begin{array}{l}
z_1^{12}+27z_2^6+9z_1^4z_2^4+8z_3^4+24z_6^2+
18z_2^2z_4^2+6z_1^4z_4^2+3z_1^8z_2^2\end{array}\right),\\[5pt]
Z_{H_5^5}({z})&=\frac{1}{120}\left(\begin{array}{l}24z_5^2+30z_2z_4^2+20z_1z_3z_6+20z_1z_3^3+
15z_1^2z_2^4+10z_1^4z_2^3+z_1^{10}\end{array}\right).
\end{align*}
Consequently, the values $F_5(k)$ for $8<k\leq 16$ can be derived from  (\ref{cg103}), and they agree with the results of Aichholzer \cite{Aic2}, see Table \ref{t21}.
\begin{table}[h,t]
\begin{center}
\begin{tabular}{|l|l|l|l|l|l|l|l|l|l|l|l|l|l|}\hline
$k$&9&10&11&12&13&14&15&16\\\hline
$H_5^1$&56&50&27&19&6&4&1&1\\\hline
$H_5^2$&159&135&68&43&12&7&1&1\\\hline
$H_5^3$&9&5&1&1&&&&\\\hline $H_5^4$&7&5&1&1&&&&\\\hline
$H_5^5$&1&1&&&&&&\\\hline
$F_5(k)$&8781&19767&37976&65600&98786&133565&158656&159110\\\hline
\end{tabular}
\end{center}
\caption{$F_5(k)$ for $8< k\leq
16$.}\label{t21}
\end{table}

The main objective of this section is to compute $F_6(k)$ for $16<k\leq 32$. As mentioned in Section 4, there are $6$ representatives of equivalence classes of spanned hyperplanes of $Q_6$ containing more than $16$ vertices  of $Q_6$, i.e., $H_6^1, H_6^2, H_6^3, H_6^4, H_6^5, H_6^6$. Again, by applying the techniques in Section 6, we obtain that
\begin{align*}
Z_{H_6^1}({z})&=Z_5({z}),\\[5pt]
Z_{H_6^2}({z})&=\frac{1}{768}\left(\begin{array}{l}
z_1^{32}+12z_1^{16}z_2^8+12z_1^8z_2^{12}+127z_2^{16}+
32z_1^8z_3^8+\\[5pt]
48z_1^4z_2^2z_4^6+168z_4^8+224z_2^4z_6^4+96z_8^4+48z_2^4z_4^6
\end{array}\right),\\[7pt]
Z_{H_6^3}({z})&=\frac{1}{288}\left(\begin{array}{l}
z_1^{24}+6z_1^{12}z_2^6+52z_2^{12}+18z_3^8+
48z_4^6+32z_2^3z_6^3+3z_1^8z_2^8+\\[5pt]
18z_1^4z_2^{10}+24z_1^2z_3^2z_2^2z_6^2+8z_1^6z_3^6+12z_3^4z_6^2+42z_6^4+
24z_{12}^2
\end{array}\right),\\[5pt]
Z_{H_6^4}({z})&=\frac{1}{384}\left(\begin{array}{l}
z_1^{24}+81z_2^{12}++2z_1^{12}z_2^6+18z_1^4z_2^{10}+
15z_1^8z_2^8+72z_6^4+32z_{12}^2\\[5pt]64z_4^6+16z_3^4z_6^2+8z_3^8
+54z_2^4z_4^4+12z_1^4z_2^2z_4^4+6z_1^8z_4^4+3z_1^{16}z_2^4
\end{array}\right),\\[7pt]
Z_{H_6^5}({z})&=\frac{1}{240}\left(\begin{array}{l}
z_1^{20}+24z_{10}^2+60z_2^2z_4^4+26z_2^{10}+20z_1^2z_3^2z_6^2+\\[5pt]
20z_1^2z_3^6+15z_1^4z_2^8+10z_1^8z_2^6+40z_2z_6^3+24z_5^4
\end{array}\right),\\[7pt]
Z_{H_6^6}({z})&=\frac{1}{1440}\left(\begin{array}{l}
z_1^{20}+144z_5^4+
144z_{10}^2+320z_2z_6^3+270z_2^2z_4^4+76z_2^{10}\\[5pt]
+90z_1^4z_4^4+30z_1^8z_2^6+45z_1^4z_2^8+240z_1^2z_3^2z_6^2+80z_1^2z_3^6
\end{array}\right).
\end{align*}
Based on relation (\ref{cg103}), we can compute $F_6(k)$ for $16<k\leq 32$. These
values are listed in Table \ref{t42}.
\begin{table}[h,t]
\centering
\begin{tabular}{|l|l|l|l|l|l|l|l|l|l|l|l|l|l|l|l|l|l|l|l|l|l|l|l|}\hline
&$H_6^1$&$H_6^2$&$H_6^3$&$H_6^4$&$H_6^5$&$H_6^6$&$F_6(k)$\\\hline
17&158658&767103&1464&1334&12&5&30063520396   \\\hline
18&133576&642880&657&630&5&3&78408664654   \\\hline
19&98804&474635&220&216&1&1&189678190615  \\\hline
20&65664&312295&81&86&1&1&426539396250  \\\hline
21&38073&179829&19&20&&&893345853436\\\hline
22&19963&92309&7&8&&&1745593621167\\\hline
23&9013&40948&1&1&&&3186944223591\\\hline
24&3779&16335&1&1&&&5443544457875\\\hline
25&1326&5500&&&&&8708686176141\\\hline
26&472&1753&&&&&13061946974320\\\hline
27&131&441&&&&&18382330104124\\\hline
28&47&129&&&&&24289841497705\\\hline
29&10&23&&&&&30151914536900\\\hline
30&5&9&&&&&35176482187384\\\hline
31&1&1&&&&&38580161986424\\\hline
32&1&1&&&&&39785643746724\\\hline
\end{tabular}
\caption{$F_6(k)$ for $16<k\leq 32$.}\label{t42}
\end{table}

\section{$H_n(k)$ for $2^{n-3}<k\leq 2^{n-2}$}

In this section, we shall present an approach for computing  $H_n(k)$ for
$2^{n-3}<k\leq 2^{n-2}$.  This enables us to
  determine $F_6(k)$ for $k=13, 14, 15, 16$.
Together with the computation of
Aichholzer up to 12 vertices for $n=6$,
we have completed the enumeration of
full-dimensional $0/1$-equivalence classes of the
6-dimensional hypercube.

Let us recall the map $\Phi$ defined in Section 4, which
will be used in the computation of $H_n(k)$  for
$2^{n-3}<k\leq 2^{n-2}$. Let $H_1,H_2,\ldots,H_{h(n,k)}$ be the representatives of equivalence classes of spanned hyperplanes of $Q_n$ containing at least $k$ vertices.  As before, denote by $\mathcal{P}(H_i,k)$ ($1\leq i\leq h(n,k)$) the set of partial $0/1$-equivalence classes of $H_i$ with $k$ vertices. Let  $\mathcal{P}_i$ be a partial $0/1$-equivalence class in $\mathcal{P}(H_i,k)$ ($1\leq i\leq h(n,k)$). So $\Phi$ maps $\mathcal{P}_i$ to the (unique) $0/1$-equivalence class in $\mathcal{H}_n(k)$ containing $\mathcal{P}_i$. When $2^{n-2}<k\leq 2^{n-1}$, it has been shown in Theorem \ref{cg16} that  $\Phi$ is  a bijection. However,
as pointed out after the proof of Theorem \ref{cg16}, when  $k\leq 2^{n-2}$, $\Phi$ is surjective but not necessarily injective.

For the purpose of computing $H_n(k)$ for
$2^{n-3}<k\leq 2^{n-2}$, we shall first derive an expression for $H_n(k)$, which is valid  for general $k$.
Let $1\leq i\leq h(n,k)$, and define \[A_i=\Phi(\mathcal{P}(H_i,k)).\] Since $\Phi$ is surjective, we see that
\[\mathcal{H}_n(k)=A_1\cap A_2\cup\cdots\cup A_{h(n,k)}.\]   It follows from the  principle of inclusion-exclusion that
\begin{equation}\label{mainformula}
\begin{split}
H_n(k)=&\sum_{1\leq i\leq h(n,k)}|A_i|-
\sum_{1\leq i_1<i_2\leq h(n,k)}|A_{i_1}\cap A_{i_2}|\\
&+\sum_{1\leq i_1<i_2<i_3 \leq h(n,k)}|A_{i_1}\cap A_{i_2}\cap A_{i_3}|-\cdots.
\end{split}
\end{equation}
Hence the task of  computing $H_n(k)$ reduces to evaluating
$|A_{i_1}\cap A_{i_2}\cap\cdots\cap A_{i_m}|$ for $1\leq i_1<i_2<\cdots<i_m\leq h(n,k)$.

Assume that $2^{n-3}<k\leq 2^{n-2}$. In what follows, we shall focus on the computation of the cardinalities of $A_i$ for $1\leq i\leq h(n,k)$, and the cardinalities of $A_i\cap A_j$ for $1\leq i<j\leq h(n,k)$. The computation for the cardinalities of  $A_{i_1}\cap A_{i_2}\cap\cdots\cap A_{i_m}$ in the general case
can be carried out in the same way.
When $n=6$ and $k=13,14,15,16$, the computations turn out to be  quite simple.

We first compute $|A_i|$ ($1\leq i\leq h(n,k)$) for $2^{n-3}<k\leq 2^{n-2}$. Since $A_i=\Phi(\mathcal{P}(H_i,k))$, we have   $|A_i|=|\mathcal{P}(H_i,k)|$. Recall that $|\mathcal{P}(H_i,k)|$
is defined as $N_{H_i}(k)$ in Section 4 and has been computed
for the case $2^{n-2}<k\leq 2^{n-1}$. To compute $N_{H_i}(k)$
for  $2^{n-3}<k\leq 2^{n-2}$,  we need some notation.

Let $H$ be a spanned hyperplane  of $Q_n$, and $\mathcal{S}$ be a subset  of $H$. Recall that $\mathcal{S}(Q_n)$ is the set of $0/1$-polytopes of $Q_n$ contained in $\mathcal{S}$.
In Section 4, we  defined the partial $0/1$-equivalence relation on $\mathcal{S}(Q_n)$. Here we need introduce another equivalence relation on $\mathcal{S}(Q_n)$, that is, two $0/1$-polytopes in $\mathcal{S}(Q_n)$ are said to be equivalent if one can be transformed to the other by a symmetry in $F(H)$. The associated equivalence classes in $\mathcal{S}(Q_n)$  are called local $0/1$-equivalence classes of $\mathcal{S}$. Since $F(H)$ is a subgroup of $B_n$,
each  local $0/1$-equivalence class of $\mathcal{S}$ is contained in a (unique) partial  $0/1$-equivalence class of $\mathcal{S}$.

Denote by $\mathcal{L}(\mathcal{S},k)$ the set of local $0/1$-equivalence classes of $\mathcal{S}$ with $k$ vertices.
When $\mathcal{S}=H$, $\mathcal{L}(H,k)$ has appeared in Section 4, that is, $\mathcal{L}(H,k)$ is the set of equivalence classes of $0/1$-polytopes contained in $H$ with $k$ vertices under the action of $F(H)$. So we have the following relation
\begin{equation}\label{local}
|\mathcal{L}(H,k)|=\left[u_1^ku_2^{|V_n(H)-k|}\right]C_H(u_1,u_2).
\end{equation}

In order to compute  $N_H(k)$ for $2^{n-3}<k\leq 2^{n-2}$, we shall define a partition of $\mathcal{L}(H,k)$ into two subsets $\mathcal{L}_\ast(H,k)$ and $\mathcal{L}^\ast(H,k)$. This
requires a property as given in Theorem \ref{L1}.

Let $H$ be a spanned hyperplane of $Q_n$ containing at least $k$ vertices of $Q_n$.
Denote by $E(H,k)$ the set of intersections  $H\cap w(H)$ such that
\begin{itemize}
\item[(1).] The symmetry $w$ of $Q_n$ does not fix $H$, that is, $H\neq w(H)$;
\item[(2).] The intersection $H\cap w(H)$ contains at least $k$ vertices of $Q_n$.
\end{itemize}
 Denote by $h_1(H,k)$ the number of equivalence classes of $E(H,k)$ under the
symmetries in $F(H)$. Let $E_1(H,k)=\{H\cap w_{i}(H)\colon 1\leq i\leq h_1(H,k)\}$ be the set of representatives of these
equivalence classes of $E(H,k)$.

Consider the (disjoint) union of  $\mathcal{L}(H\cap w_i(H),k)$, where $1\leq i\leq h_1(H,k)$. We shall define a map $\Phi_1$ from this union to $\mathcal{L}(H,k)$.
For $1\leq i\leq h_1(H,k)$, let  $\mathcal{L}_i$ be a local $0/1$-equivalence class in $\mathcal{L}(H\cap w_i(H),k)$. Evidently,  there is a (unique) local $0/1$-equivalence class in $\mathcal{L}(H,k)$ containing  $\mathcal{L}_i$, denoted $\mathcal{L}_i'$. Define $\Phi_1(\mathcal{L}_i)=\mathcal{L}_i'$. Then we have the following property.

\begin{theorem}\label{L1}
If $2^{n-3}<k\leq 2^{n-2}$, then the map $\Phi_1$ is an injection.
\end{theorem}

\pf   Let $\mathcal{L}$ and $\mathcal{L}'$ be two distinct local $0/1$-equivalence classes with $k$ vertices. Assume that $\mathcal{L}$ (resp. $\mathcal{L}'$) is in $\mathcal{L}(H\cap w_i(H),k)$ (resp. $\mathcal{L}(H\cap w_j(H),k)$), where $1\leq i, j\leq h_1(H,k)$. To prove that $\Phi_1$ is an injection, we need to show that $\Phi_1(\mathcal{L})\neq \Phi_1(\mathcal{L}')$. Clearly, if $i=j$ then we see that  $\Phi_1(\mathcal{L})\neq \Phi_1(\mathcal{L}')$. We now consider the case $i\neq j$.

Assume to the contrary that $\Phi_1(\mathcal{L})= \Phi_1(\mathcal{L}')$.
Let $P$ (resp. $P'$) be any $0/1$-polytope in $\mathcal{L}$ (resp. $\mathcal{L}'$). Then there is a symmetry $w\in F(H)$ such that $P=w(P')$. Since both $P$ and $P'$  have more than $2^{n-3}$ vertices of $Q_n$, we see from Theorem \ref{cg14} that $\mathrm{dim}(P)=\mathrm{dim}(P')\geq n-2$.
Since $P$ (resp. $P'$) is contained in $H\cap w_i(H)$ (resp. $H\cap w_j(H)$),  both $P$ and $P'$ are of dimension  $n-2$. This implies that $H\cap w_i(H)$ (resp. $H\cap w_j(H)$) is the affine space spanned by  $P$ (resp. $P'$). Hence we deduce that $H\cap w_i(H)=w(H\cap w_j(H))$, which is contrary to the assumption that $H\cap w_i(H)$ and $H\cap w_j(H)$ are not equivalent under the symmetries in $F(H)$. This completes the proof. \qed

We are now ready to  define $\mathcal{L}_\ast(H,k)$
to be the image of $\Phi_1$. More precisely, $\mathcal{L}_\ast(H,k)$ is the  (disjoint)  union of $\Phi_1(\mathcal{L}(H\cap w_i(H),k))$, where $1\leq i\leq h_1(H,k)$. Let
\begin{equation}\label{E1}
\mathcal{L}^{\ast}(H,k)=\mathcal{L}(H,k)\backslash \mathcal{L}_\ast(H,k).\end{equation}
From the above definition (\ref{E1}), it can be seen that, for any local $0/1$-equivalence class $\mathcal{L}\in \mathcal{L}^{\ast}(H,k)$ and any $0/1$-polytope  $P\in \mathcal{L}$, if $w\in B_n$ is a symmetry such that $w(P)$ is contained in $H$, then $w(H)=H$. This yields that
$\mathcal{L}$ is also a partial $0/1$-equivalence class of $H$.
Consequently,  $\mathcal{L}^\ast(H,k)$ is a subset of $\mathcal{P}(H,k)$.
Let \begin{equation}\label{E2}
\mathcal{P}_{\ast}(H,k)=\mathcal{P}(H,k)\backslash\mathcal{L}^\ast(H,k).
\end{equation}
Combining (\ref{local}), (\ref{E1}) and (\ref{E2}),  we find that
\begin{equation}\label{identity}
\begin{split}
N_H(k)&=|\mathcal{P}(H,k)|\\
&=|\mathcal{L}^\ast(H,k)|+|\mathcal{P}_{\ast}(H,k)|\\[5pt]
&=|\mathcal{L}(H,k)|-|\mathcal{L}_\ast(H,k)|+|\mathcal{P}_{\ast}(H,k)|\\[5pt]
&=\left[u_1^ku_2^{|V_n(H)-k|}\right]C_H(u_1,u_2)-|\mathcal{L}_\ast(H,k)|+
|\mathcal{P}_{\ast}(H,k)|.
\end{split}
\end{equation}
Therefore, for $2^{n-3}<k\leq 2^{n-2}$  $N_H(k)$ is determined by  the cardinalities of $\mathcal{L}_\ast(H,k)$ and $\mathcal{P}_{\ast}(H,k)$.
From Theorem \ref{L1}, we see that  for $2^{n-3}<k\leq 2^{n-2}$, $|\mathcal{L}_\ast(H,k)|$ can be derived from   the cardinalities of $\mathcal{L}(H\cap w(H),k)$, where $H\cap w(H)\in E_1(H,k)$. We shall demonstrate that the computation of $|\mathcal{P}_{\ast}(H,k)|$ for $2^{n-3}<k\leq 2^{n-2}$ can be  carried out in a similar fashion.

Denote by $h_2(H,k)$ the number of equivalence classes of $E(H,k)$ under the symmetries of $Q_n$.  Let \[ E_2(H,k)=\{H\cap w_{i}(H)\colon 1\leq i\leq h_2(H,k)\}\] be the set of representatives of these equivalence classes of $E(H,k)$.
We define a map $\Phi_2$ from  the (disjoint) union of $\mathcal{P}(H\cap w_i(H),k)$, where $1\leq i\leq h_2(H,k)$,
 to $\mathcal{P}_{\ast}(H,k)$.
Let  $\mathcal{P}$ be a partial $0/1$-equivalence class of $\mathcal{P}(H\cap w_i(H),k)$ ($1\leq i\leq h_2(H,k)$). Then the image $\Phi_2(\mathcal{P})$ is defined to be the  (unique) partial $0/1$-equivalence class of $\mathcal{P}_\ast(H,k)$ that contains  $\mathcal{P}$. We reach the following assertion. The proof is
similar to that of Theorem \ref{L1}, hence it is omitted.

\begin{theorem}\label{L2}
If $2^{n-3}<k\leq 2^{n-2}$, then the map $\Phi_2$ is a bijection.
\end{theorem}

So far, we see that the number $N_H(k)$ for $2^{n-3}<k\leq2^{n-2}$ can be computed based on the cardinalities of $\mathcal{L}(H\cap w(H),k)$ and $\mathcal{P}(H\cap w(H),k)$, where $H\cap w(H)\in E(H,k)$.
We shall illustrate  how to compute $|\mathcal{L}(H\cap w(H),k)|$ and $|\mathcal{P}(H\cap w(H),k)|$ for $2^{n-3}<k\leq2^{n-2}$.

Assume that  $H\cap w(H)\in E(H,k)$.
Let $P$ and $P'$ be any two $0/1$-polytopes
belonging to the same local  (resp. partial)  $0/1$-equivalence class of $H\cap w(H)$ with $k$ vertices.  Then
there exists a symmetry in $F(H)$ (resp. $B_n$) such that
$w(P)=P'$. It is clear from  Theorem \ref{cg14} that
both $P$ and $P'$ have  dimension
$n-2$. Hence  $H\cap w(H)$ is the affine space spanned by $P$,
 or, equivalently, by $P'$. So we deduce that
$w(H\cap w(H))=H\cap w(H)$. This implies that for $2^{n-3}<k\leq 2^{n-2}$, we can use  P\'olya's theorem to compute the number of local (resp. partial) $0/1$-equivalence classes of $H\cap
w(H)$ with $k$ vertices.

Let
\[F(H,w)=\left\{w'\in F(H)\colon\, w'\left(H\cap w(H)\right)=H\cap w(H)\right\}\] and
\[F\left(H\cap
w(H) \right)=\left\{w'\in B_n\colon\, w'\left(H\cap
w(H)\right)=H\cap
w(H)\right\}.\]
Denote by $V_n(H\cap w(H))$  the set of vertices of $Q_n$ contained in $H\cap w(H)$, and denote by $Z_{(H,w)}({z})$ (resp. $Z_{H\cap
w(H)}({z})$) the  cycle index of $F\left(H,w\right)$  (resp. $F\left(H\cap
w(H) \right)$) acting on $V_n(H\cap w(H))$.  Write $C_{(H,w)}(u_1,u_2)$ (resp. $C_{H\cap w(H)}(u_1,u_2)$) for the polynomial obtained from $Z_{(H,w)}({z})$ (resp. $Z_{H\cap w(H)}(z)$) by substituting $z_i$ with $u_1^i+u_2^i$.
Thus, for $2^{n-3}<k\leq 2^{n-2}$, we obtain that
\begin{equation}\label{certain2}
|\mathcal{L}(H\cap w(H),k)|=\left[u_1^{k}u_2^{|V_n(H\cap w(H))|-k}\right]C_{(H,w)}(u_1,u_2)\end{equation}
and
\begin{equation}\label{certain}
|\mathcal{P}(H\cap w(H),k)|=\left[u_1^{k}u_2^{|V_n(H\cap w(H))|-k}\right]C_{H\cap w(H)}(u_1,u_2).\end{equation}
Thus, applying   and Theorems \ref{L1} and \ref{L2} and plugging the above   formulas (\ref{certain2}) and (\ref{certain}) into (\ref{identity}), we arrive at the following relation.

\begin{theorem}\label{cg102}
Let $2^{n-3}<k\leq 2^{n-2}$,and
$H$ be a spanned hyperplane of $Q_n$ containing at least $k$ vertices of $Q_n$.  Set $q(w)=|V_n(H\cap w(H))|$. Then we have
\begin{equation}\label{formula}
\begin{split}
N_H(k)=&\left[u_1^{k}u_2^{|V_n(H)|-k}\right]C_H(u_1,u_2)-\sum_{H\cap w(H)\in E_1(H,k)}\left[u_1^{k}u_2^{q(w)-k}\right]C_{(H,w)}(u_1,u_2)\\[5pt]
&+\sum_{H\cap w(H)\in E_2(H,k)}\left[u_1^{k}u_2^{q(w)-k}\right]C_{H\cap w(H)}(u_1,u_2).
\end{split}\end{equation}
\end{theorem}

Theorem \ref{cg102} enables us to compute $N_H(k)$ for $k=13,14,15,16$, where $H$ is a spanned hyperplane of $Q_6$ containing more than $12$ vertices of $Q_6$.
In addition to $H_6^1, H_6^2, H_6^3, H_6^4, H_6^5, H_6^6$, we have $8$ representatives of equivalence classes of spanned hyperplanes of $Q_6$ containing more than $12$ vertices of $Q_6$, namely,
\begin{align*}
&H_1\colon x_1+x_2+x_3+2x_4=2,\\[3pt]
&H_2\colon x_1+x_2+x_3+x_4=1,\\[3pt]
&H_3\colon x_1+x_2+x_3+x_4+2x_5=3,\\[3pt]
&H_{4}\colon x_1+x_2+x_3+x_4+x_5+2x_6=3,\\[3pt]
&H_5\colon x_1+x_2+x_3+x_4+x_5+x_6=2,\\[3pt]
&H_6\colon x_1+x_2+x_3+x_4+2x_5=2,\\[3pt]
&H_7\colon x_1+x_2+x_3+2x_4+2x_5=3,\\[3pt]
&H_8\colon x_1+x_2+x_3+x_4+2x_5+2x_6=4.
\end{align*}
It is easily checked that $E(H,k)=\emptyset$ for  $k=13,14,15,16$, except for the  two spanned hyperplanes  $H_6^1$ and $H_6^2$. Therefore, we  can deduce from Theorem \ref{cg102}  that
\begin{equation}\label{c1}
N_H(k)=\left[u_1^{k}u_2^{V_n(H)-k}\right]C_H(u_1,u_2),
\end{equation}where $H=H_6^3-H_6^6, H_1-H_8$.  The cycle indices for $H=H_6^3-H_6^6$ have been given  in Section 7. For $H=H_6$, $H_7$ and $H_8$, it is easily verified that $N_{H}(13)=2$ and $N_{H}(14)=1$.  Using the techniques in Section 6, we can
 derive the cycle indices for $H_1-H_5$ as shown below:
\begin{align*}
Z_{H_1}({z})&=\frac{1}{48}\left(\begin{array}{l}
z_1^{16}+4z_{12}z_4+4z_3^2z_6z_1^2z_2+2z_3^4z_1^4+\\[5pt]
12z_2^8+8z_4^4+6z_1^4z_2^6+5z_1^8z_2^4+6z_6^2z_2^2
\end{array}\right),\\[7pt]
Z_{H_2}({z})&=\frac{1}{192}\left(\begin{array}{l}
z_1^{16}+68z_4^4+24z_6^2z_2^2+16z_{12}z_4
+8z_3^4z_1^4\\[5pt]+39z_2^8+12z_1^4z_2^6+8z_1^8z_2^4+16z_3^2z_6z_1^2z_2
\end{array}\right),\\[5pt]
Z_{H_3}({z})&=\frac{1}{96}\left(\begin{array}{l}
z_1^{16}+24z_6^2z_2^2+8z_3^4z_1^4+33z_2^8+ 6z_1^8z_2^4+24z_4^4
\end{array}\right),\\[7pt]
Z_{H_4}({z})&=\frac{1}{120}\left(\begin{array}{l}
z_1^{15}+24z_5^3+30z_2z_4^3z_1+20z_1z_3^2z_6z_2+
20z_1^3z_3^4+15z_1^3z_2^6+10z_1^7z_2^4
\end{array}\right),\\[7pt]
Z_{H_5}({z})&=\frac{1}{720}\left(\begin{array}{l}
z_1^{15}+120z_3z_6^2+144z_5^3+40z_3^5+180z_1z_2z_4^3\\[5pt]+40z_1^3z_3^4+60z_1^3z_2^6+
15z_1^7z_2^4+ 120z_1z_2z_3^2z_6
\end{array}\right).
\end{align*}

It remains to compute $N_H(k)$ for $H=H_6^1$ and $H_6^2$ for $k=13,14,15,16$.
For $H_6^1\colon x_1=0$ and $k=13,14,15,16$, it is routine to check that
\begin{align*}
E_1(H_6^1,k)=E_2(H_6^1,k)
=\left\{H_6^1\cap w(H_6^1)\colon w=(1,2)(3)(4)(5)(6)\right\},
\end{align*} that is,
\begin{equation*}
\{(x_1,\ldots,x_6)\colon x_1=0, x_2=0\}.\end{equation*}
Thus, for $k=13,14,15,16$, it is clear that both the numbers of  local and partial $0/1$-equivalence classes of $H_6^1\cap w(H_6^1)$ with $k$ vertices are given by
\[\left[u_1^{k}u_2^{16-k}\right]C_4(u_1,u_2).\] Therefore, for $k=13,14,15,16$,  by Theorem \ref{cg102} we find that
\begin{equation}\label{r1}
N_{H_6^1}(k)=\left[u_1^{k}u_2^{32-k}\right]C_{H_6^1}(u_1,u_2).
\end{equation}

Finally, we come to the computation of $N_{H_6^2}(k)$ for $k=13,14,15,16$. In this case, it is easy to check  that
\begin{equation*}
E_1(H_6^2,k)=E_2(H_6^2,k)=\{H_6^2\cap w_1(H_6^2),H_6^2\cap w_2(H_6^2)\},
\end{equation*}
where $w_1=(1,3,2)(4)(5)(6)$ and $w_2=(1,3)(2,4)(5)(6)$.
Since
\begin{align*}
V_6(H_6^2\cap w_1(H_6^2))=&\{(1,0,1,x_4,x_5,x_6)\colon \mbox{$x_i=0$ or 1 for $i=4,5,6$}\}\cup \\
&\{(0,1,0,x_4,x_5,x_6)\colon \mbox{$x_i=0$ or 1 for $i=4,5,6$}\},
\end{align*}
it can be easily checked  that $\mathcal{L}(H_6^2\cap w_1(H_6^2),k)=\mathcal{P}(H_6^2\cap w_1(H_6^2),k)$ for $k=13,14,15,16$. By Theorem \ref{cg102}, we obtain that for
$k=13,14,15,16$,
\begin{equation}\label{r2}
\begin{split}
N_{H_6^2}=&\left[u_1^{k}u_2^{32-k}\right]C_{H_6^2}(u_1,u_2)-
\left[u_1^{k}u_2^{16-k}\right]C_{(H_6^2,w_2)}(u_1,u_2)\\[5pt]
&+
\left[u_1^{k}u_2^{16-k}\right]C_{H_6^2\cap w_2(H_6^2)}(u_1,u_2).
\end{split}
\end{equation}

Next, we proceed to demonstrate how to compute $|A_i\cap A_j|$ for $1\leq i<j\leq h(n,k)$. Let $E(H_i,H_j,k)$ be the set of intersections  $H_i\cap w(H_j)$ ($w\in B_n$) that contain at least $k$ vertices of $Q_n$.
Denote by $h(H_i,H_j,k)$ the number of equivalence classes of $E(H_i,H_j,k)$ under the symmetries of $Q_n$.
Let $m=h(H_i,H_j,k)$. Assume that $E_1(H_i,H_j)=\{H_i\cap w_1(H_j),\ldots,H_i\cap w_m(H_j)\}$ is the set of  representatives of equivalence classes in $E(H_i,H_j,k)$. We define a map $\Phi_3$ from the union of $\mathcal{P}(H_i\cap w_s(H_j),k)$, where $1\leq s\leq h(H_i,H_j,k)$, to $A_i\cap A_j$. Let  $\mathcal{P}_s$ be a partial $0/1$-equivalence class in $\mathcal{P}(H_i\cap w_s(H_j),k)$. Clearly,  there is a (unique) partial $0/1$-equivalence class in $A_i\cap A_j$ containing  $\mathcal{P}_s$, which will be denoted by $\mathcal{P}_s'$. Define $\Phi_3(\mathcal{P}_s)=\mathcal{P}_s'$. We have the following conclusion. We omit the proof since it is
similar to that of Theorem \ref{L1}.

\begin{theorem}\label{L3}
If $2^{n-3}<k\leq 2^{n-2}$, then the map $\Phi_3$ is a bijection.
\end{theorem}

As a  consequence of Theorem \ref{L3}, for $2^{n-3}<k\leq 2^{n-2}$, we have
\begin{equation*}
|A_i\cap A_j|=\sum_{s=1}^{h(H_i,H_j,k)}|\mathcal{P}(H_i\cap w_s(H_j),k)|.
\end{equation*}

The computation for $|A_{i_1}\cap A_{i_2}\cap\cdots\cap A_{i_m}|$ ($m\geq 3$) in the general case can be done in a similar fashion. In fact, it will be shown that for $2^{n-3}<k\leq 2^{n-2}$, the computation  can be reduced to the case $m=2$.

Let $2^{n-3}<k\leq 2^{n-2}$, and $E(H_{i_1},\ldots,H_{i_m},k)$ be the set of intersections  $H_{i_1}\cap w_2(H_{i_2})\cap \cdots \cap w_m(H_{i_m})$, where $w_i$ for $2\leq i\leq m$ are symmetries of $Q_n$,  that contain at least $k$ vertices of $Q_n$. Denote by $E_1(H_{i_1},\ldots,H_{i_m},k)$ the set of representatives of equivalence classes of $E(H_{i_1},\ldots,H_{i_m},k)$ under the symmetries of $Q_n$. We define a map $\Phi_m$ from the (disjoint) union of $\mathcal{P}(H_{i_1}\cap w_2(H_{i_2})\cap \cdots \cap w_m(H_{i_m}),k)$, where \[H_{i_1}\cap w_2(H_{i_2})\cap \cdots \cap w_m(H_{i_m})\in E_1(H_{i_1},\ldots,H_{i_m},k),\] to the set $A_{i_1}\cap A_{i_2}\cap\cdots\cap A_{i_m}$.
Let $\mathcal{P}\in \mathcal{P}(H_{i_1}\cap w_2(H_{i_2})\cap \cdots \cap w_m(H_{i_m}),k)$. The image $\Phi_m(\mathcal{P})$ is defined to be the unique partial $0/1$-equivalence class in $A_{i_1}\cap A_{i_2}\cap\cdots\cap A_{i_m}$ containing $\mathcal{P}$.
Similarly,  we can prove that if $2^{n-3}<k\leq 2^{n-2}$, then $\Phi_m$ is  a bijection. Thus we deduce that for  $2^{n-3}<k\leq 2^{n-2}$,
\[|A_{i_1}\cap A_{i_2}\cap\cdots\cap A_{i_m}|=\sum|\mathcal{P}(H_{i_1}\cap w_2(H_{i_2})\cap \cdots \cap w_m(H_{i_m}),k)|,\] where the sum ranges over the representatives of $E_1(H_{i_1},\ldots,H_{i_m},k)$.

We further claim  for  $2^{n-3}<k\leq 2^{n-2}$,
$E(H_{i_1},\ldots,H_{i_m},k)$ is a subset of $E(H_{i_1},E_{i_2})$.
This can be proved as follows. Assume that
$2^{n-3}<k\leq 2^{n-2}$, and that $H_{i_1}\cap w_2(H_{i_2})\cap
\cdots \cap w_m(H_{i_m})$ is in $E(H_{i_1},\ldots,H_{i_m},k)$.
From Theorem \ref{cg14} it can be seen that the dimension of
$H_{i_1}\cap w_2(H_{i_2})\cap \cdots \cap w_m(H_{i_m})$ is at
least  $n-2$, since it contains more than $2^{n-3}$ vertices of
$Q_n$. On the other hand, it is clear that $H_{i_1}\cap
w_2(H_{i_2})\cap \cdots \cap w_m(H_{i_m})$ has dimension at most
$n-2$. Hence, when  $2^{n-3}<k\leq 2^{n-2}$, we conclude that
$H_{i_1}\cap w_2(H_{i_2})\cap \cdots \cap w_m(H_{i_m})$ is of
dimension $n-2$. Hence we obtain that $H_{i_1}\cap w_2(H_{i_2})\cap
\cdots \cap w_m(H_{i_m})=H_{i_1}\cap w_2(H_{i_2})$. Therefore,
$E(H_{i_1},\ldots,H_{i_m},k)$ is a subset of $E(H_{i_1},E_{i_2})$.
This implies that the for $2^{n-3}<k\leq 2^{n-2}$, the computation
for $E(H_{i_1},\ldots,H_{i_m},k)$ can be reduced
to the case $m=2$.
More specifically, for  $2^{n-3}<k\leq 2^{n-2}$, an
intersection $H_{i_1}\cap w_2(H_{i_2})\cap \cdots \cap
w_m(H_{i_m})$ belongs to $E(H_{i_1}, \cdots , H_{i_m},k)$
whenever (possibly after the action of some symmetry of $Q_n$) it
belongs to  $E(H_{i_{j_1}},H_{i_{j_2}})$ for $1\leq
{j_1}<{j_2}\leq m$.

We now turn to the case when $n=6$ and  $k=13,14,15,16$. All possible pairs $\{H_i,H_j\}$  such that $E(H_i,H_j,k)$ is nonempty are listed below.

(1). $\{H_6^1, H_6^2\}$. In this case, it can be easily checked that
\begin{equation}
\begin{split}
E_1(H_6^1, H_6^2,k)&=\{H_6^1\cap H_6^2\}\cup \{H_6^1\cap w(H_6^2)\colon w=(1,3,2)(4)(5)(6)\}\\
&=\{H_6^1\cap H_6^2\}\cup \{H_6^1\cap H_6^3\}.
\end{split}
\end{equation}

(2). $\{H_6^1,H_6^3\}$ and $\{H_6^2,H_6^3\}$.
In these two cases, we have
\begin{equation}
E_1(H_6^1, H_6^3,k)=E_1(H_6^2\cap H_6^3,k)=\{H_6^1\cap H_6^3\}.
\end{equation}

(3). $\{H_6^2, H_6^4\}$. In this case, it can be verified that
\begin{equation}
E_1(H_6^2, H_6^4,k)=\{H_6^2\cap H_6^4\}.
\end{equation}

From the above, we see that $H_6^1$, $H_6^2$ and $H_6^3$ are the
only hyperplanes such that for $k=13,14,15,16$,
$E(H_{i_1},H_{i_2},H_{i_3},k)$ is nonempty. Moreover,  for
$k=13,14,15,16$ we have
\begin{equation}
E_1(H_6^1,H_6^2,H_6^3,k)=\{H_6^1\cap H_6^3\}.
\end{equation}

For $k=13,14,15,16$, it is easy to see that
\begin{equation}\label{s1}
\begin{split}
|\mathcal{P}(H_6^1\cap H_6^2,k)|&=\left[u_1^{k}u_2^{16-k}\right]C_4(u_1,u_2),\\
|\mathcal{P}(H_6^1\cap H_6^3,k)|&=\left[u_1^{k}u_2^{16-k}\right]C_{H_5^2}(u_1,u_2),\\
|\mathcal{P}(H_6^2\cap H_6^4,k)|&=\left[u_1^{k}u_2^{16-k}\right]C_{H_6^2\cap w(H_6^2)}(u_1,u_2),
\end{split}
\end{equation}where $w=(1,3)(2,4)(5)(6)$.

From (\ref{mainformula}) and the relations (\ref{c1})--(\ref{s1}), we deduce that for $n=6$ and $k=13,14,15,16$,
\begin{equation}\label{rr}
\begin{split}
H_6(k)=&\sum_{i=1}^6\left[u_1^{k}u_2^{|V_6(H_6^i)|-k}\right]C_{H_6^i}(u_1,u_2)
+\sum_{i=1}^8\left[u_1^{k}u_2^{|V_6(H_i)|-k}\right]C_{H_i}(u_1,u_2)\\
&-\left[u_1^{k}u_2^{16-k}\right]C_4(u_1,u_2)-
2\left[u_1^{k}u_2^{16-k}\right]C_{H_5^2}(u_1,u_2)-
\left[u_1^{k}u_2^{16-k}\right]C_{(H_6^2,w)}
\end{split},
\end{equation}where $w=(1,3)(2,4)(5)(6)$.
Using the argument in Section 6, for $w=(1,3)(2,4)(5)(6)$ we obtain that
\begin{equation}\label{s5}
Z_{(H_6^2,w)}=\frac{1}{32}
\left(z_1^{16}+21z_2^8+8z_4^4+2z_1^8z_2^4\right).
\end{equation}
Hence, from (\ref{rr}) and (\ref{s5}) we obtain the values of   $H_6(k)$ for $k=13,14,15,16$. Utilizing the relation $F_6(k)=A_6(k)-H_6(k)$, we deduce $F_6(k)$ for $k=13,14,15,16$ as given  in Table \ref{t5}.
\begin{table}[h,t]
\begin{center}
\begin{tabular}{|l|l|l|l|l|l|}\hline
$k$&13&14&15&16\\\hline
$F_6(k)$&290159817&1051410747&3491461629&10665920350\\\hline
\end{tabular}
\end{center}
\caption{$F_6(k)$ for $k=13,14,15,16$.}\label{t5}
\end{table}

\vspace{.2cm} \noindent{\bf Acknowledgments.} This work was
supported by  the 973 Project, the PCSIRT Project of the Ministry
of Education, and the National Science Foundation of China.

\end{document}